\documentclass[11pt]{article}
\usepackage[a4paper,margin=1in]{geometry}
\usepackage[utf8]{inputenc} 
\usepackage{lineno}
\usepackage{hyperref}
\modulolinenumbers[5]

\usepackage[T1]{fontenc}
\usepackage{amsmath,color}
\usepackage{amssymb}  
\usepackage{amsfonts}
\usepackage{mathtools}
\usepackage{graphicx}
\usepackage{dsfont}
\usepackage{psfrag}
\usepackage{subfigure,epstopdf}
\usepackage{bbm}
\usepackage{array}

\usepackage{algorithmicx}
\usepackage{algorithm}
\usepackage{algpseudocode}
\usepackage{float}
\usepackage{wrapfig}
\usepackage{enumerate}

\usepackage{xcolor}

\usepackage{booktabs}   

\usepackage{bm}

\usepackage{natbib}

\usepackage{slashed}
\usepackage{tikz,gincltex}
\usepackage{arydshln}
\usepackage{nccmath}
\usepackage{caption}
\usepackage{MnSymbol}
\usepackage{pgfplots,filecontents}
\pgfplotsset{compat=newest}
\usetikzlibrary{positioning,chains,fit,shapes,calc}
\usetikzlibrary{backgrounds}
\usetikzlibrary{shapes,snakes}
\usetikzlibrary{shapes.geometric, arrows.meta}
\pgfplotsset{compat=1.18}
\usepgfplotslibrary{groupplots}

\tikzstyle{arrow} = [thick,->,>=stealth]

\usepgfplotslibrary{statistics}
\usetikzlibrary{decorations.pathmorphing}

\makeatletter
\renewcommand*{\@begintheorem}[2]{\trivlist
      \item[\hskip \labelsep{\bfseries #1\ #2}]\itshape}
\renewcommand*{\@opargbegintheorem}[3]{\trivlist
      \item[\hskip \labelsep{\bfseries #1\ #2}] {(#3)}\ \itshape}
\makeatother

\definecolor{c1}{rgb}{0.368417, 0.506779,0.709798}
\definecolor{c2}{rgb}{0.880722, 0.611041,0.142051}
\definecolor{c3}{rgb}{0.560181, 0.691569,0.194885}
\definecolor{c4}{rgb}{0.922526, 0.385626,0.209179}
\definecolor{c5}{rgb}{0.528488, 0.470624,0.701351}
\definecolor{c6}{rgb}{0.772079, 0.431554,0.102387}
\definecolor{c7}{rgb}{0.363898, 0.618501,0.782349}
\definecolor{ao(english)}{rgb}{0.0, 0.5, 0.0}
\definecolor{cadmiumred}{rgb}{0.89, 0.0, 0.13}
\definecolor{forestgreen}{rgb}{0.13, 0.55, 0.13}
\definecolor{applegreen}{rgb}{0.55, 0.71, 0.0}

\newcommand{\false}{\textsc{false}}
\newcommand{\true}{\textsc{true}}

\tikzstyle{startstop} = [rectangle, rounded corners, minimum width=3cm, minimum height=1cm, text centered, draw=black, fill=c1]
\tikzstyle{process} = [rectangle, minimum width=3cm, minimum height=1cm, text centered, draw=c2, fill=c2!40]
\tikzstyle{io} = [rectangle, rounded corners, minimum width=3cm, minimum height=1cm, text centered, draw=c3, fill=c3!40]
\tikzstyle{io2} = [rectangle, rounded corners, minimum width=3cm, minimum height=1cm, text centered, draw=c3, fill=c3!40]
\tikzstyle{io3} = [diamond, rounded corners, minimum width=3cm, minimum height=1cm, text centered, draw=c4, fill=c4!40]
\tikzstyle{io4} = [rectangle, minimum width=3cm, minimum height=1cm, text centered, draw=c1, fill=c1!40]

\newfloat{alg}{htbp}{loa}
\floatname{alg}{Algorithm}

\newenvironment{myproof}{%
    \par\medskip\noindent%
    \textit{Proof.}\enspace\ignorespaces%
}{%
    \hfill$\square$\par\medskip%
}

\newtheorem{theorem}{Theorem}
\newtheorem{lemma}{Lemma}
\newtheorem{proposition}{Proposition}
\newtheorem{definition}{Definition}
\newtheorem{corollary}{Corollary}
\newtheorem{remark}{Remark}

\definecolor{ao}{rgb}{0.55, 0.71, 0.0}
\definecolor{bleudefrance}{rgb}{0.19, 0.55, 0.91}
\definecolor{dimgray}{rgb}{0.41, 0.41, 0.41}    
\definecolor{mediumorchid}{rgb}{0.73, 0.33, 0.83}
\definecolor{mediumtealblue}{rgb}{0.0, 0.33, 0.71}
\definecolor{harvestgold}{rgb}{0.85, 0.57, 0.0}
\definecolor{blue(pigment)}{rgb}{0.2, 0.2, 0.6}
\definecolor{forestgreen(traditional)}{rgb}{0.27, 0.35, 0.27}
\definecolor{cadmiumred}{rgb}{0.89, 0.0, 0.13}
\definecolor{orange(webcolor)}{rgb}{1.0, 0.5, 0.0}
\definecolor{tangerine}{rgb}{0.95, 0.52, 0.0}
\definecolor{c1}{rgb}{0.368417, 0.506779,0.709798}
\definecolor{c2}{rgb}{0.880722, 0.611041,0.142051}
\definecolor{c3}{rgb}{0.560181, 0.691569,0.194885}
\definecolor{c4}{rgb}{0.922526, 0.385626,0.209179}
\definecolor{c5}{rgb}{0.528488, 0.470624,0.701351}
\definecolor{c6}{rgb}{0.772079, 0.431554,0.102387}
\definecolor{c7}{rgb}{0.363898, 0.618501,0.782349}
\definecolor{c8}{rgb}{1, 0.75, 0}
\definecolor{c9}{rgb}{0.647624, 0.37816, 0.614037}
\definecolor{c10}{rgb}{0.571589, 0.586483, 0.}
\definecolor{c11}{rgb}{0.915, 0.3325, 0.2125}
\definecolor{c12}{rgb}{0.83, 0.46, 0.}
\definecolor{c13}{rgb}{0.9575, 0.545, 0.11475}
\definecolor{c14}{rgb}{1., 0.7575, 0.}
\definecolor{c15}{rgb}{0.6175, 0.715, 0.}
\definecolor{c16}{rgb}{0.15, 0.715, 0.595}
\definecolor{c17}{rgb}{0.3625, 0.545, 0.85}
\definecolor{c18}{rgb}{0.575, 0.4175, 0.85}
\definecolor{c19}{rgb}{0.677, 0.358, 0.595}
\definecolor{c20}{rgb}{0.7875, 0.358, 0.425}
\definecolor{c21}{rgb}{0.915, 0.3325, 0.2125}

\providecommand{\R}{\mathbb{R}}

\providecommand{\SG}{\Sigma_G}

\providecommand{\tr}{\operatorname{tr}}
\begin{document}

\title{Graphical stability conditions for sparse matrix spaces: a two-sided gap%
\thanks{This work is funded by national funds through FCT – Fundação para a Ciência e a Tecnologia, I.P., and, when eligible, co-funded by EU funds under project/support UID/50008/2025 – Instituto de Telecomunicações, with DOI identifier \url{https://doi.org/10.54499/UID/50008/2025}.}}

\author{%
  Diogo Po\c{c}as\thanks{Instituto Superior Técnico, Universidade de Lisboa, 1049-001 Lisbon, Portugal; and Instituto de Telecomunicações, 1049-001 Lisboa, Portugal. \texttt{diogo.pocas@tecnico.ulisboa.pt}}
  \and Guilherme Ramos\thanks{Instituto Superior Técnico, Universidade de Lisboa, 1049-001 Lisbon, Portugal; and Instituto de Telecomunicações, 1049-001 Lisboa, Portugal. \texttt{guilherme.ramos@tecnico.ulisboa.pt}}
  \and Sérgio Pequito\thanks{Instituto Superior Técnico, Universidade de Lisboa, 1049-001 Lisbon, Portugal; and Institute for Systems and Robotics, Instituto Superior Técnico, Universidade de Lisboa, Portugal. \texttt{sergio.pequito@tecnico.ulisboa.pt}}%
}
\date{}

\maketitle

\begin{abstract}
A sparse matrix space is a vector space of real matrices with a
prescribed zero pattern. It is stable if it contains a Hurwitz matrix.
Prior work by Belabbas bracketed stability between a necessary and a
sufficient graphical condition. We show that this bracket is loose on both
sides, and does not vanish even under strong connectivity: the sufficient condition
already fails to be necessary at order three, and the necessary one already
fails to be sufficient at order four. 
We propose a stronger necessary condition, using Routh-Hurwitz tables, which
catches the second example but remains strictly weaker than stability already at order five. 
We give a randomized polynomial-time algorithm for each necessary condition, via
an exact reduction to polynomial identity testing, placing them in $\mathrm{RP}$. In contrast, we prove that the sufficient condition is $\mathrm{NP}$-complete, by reduction from satisfiability, and implement a backtracking-with-memoization exponential algorithm.
We perform comprehensive numerical experiments, confirming that the gap between the necessary and sufficient conditions persists as the dimension increases, and that the algorithms retain their contrastive worst-case computational cost in practice.
\end{abstract}


\vspace{1em}

\section{Introduction}\label{sec:intro}

Deciding whether a prescribed interconnection topology can sustain stable
linear dynamics is a basic question in decentralized and networked control.
Structural systems theory addresses such questions by studying properties that
hold generically from the interconnection pattern alone, independently of the
numerical parameters~\citep{lin1974structural,dion2003generic,ramos2022overview}. 
For
controllability, the foundational results of~\citep{lin1974structural} and the
survey of~\citep{dion2003generic} characterize the property through paths and
matchings in the pattern digraph, and this line has since expanded into a broad
body of work on large-scale and networked systems~\citep{liu2011controllability,pasqualetti2014controllability,ramos2022overview},
including graphical and combinatorial methods for input selection~\citep{pequito2017robust,ramos2024minimum}.

Structural \emph{stability}, whether a pattern admits a stable realization, is comparatively underdeveloped, and appears to be less tractable to graphical analysis than controllability.
Whereas controllability of a pattern is a genericity statement, decided by a
single combinatorial witness, stability asks for a real parameter assignment
placing all eigenvalues in the open left half plane, a feasibility question over the reals rather than a purely combinatorial one. 
This shift from a combinatorial certificate to a semialgebraic feasibility problem makes the graphical characterization of stability delicate. 
Graph-theoretic
characterizations have been obtained in special settings, for instance for
interconnected homogeneous subsystems via the spectra of signed adjacency
matrices~\citep{pirani2014stability}, but a characterization for general
patterns has remained elusive.

Following~\citet{belabbas2013algorithms,belabbas2013sparse}, we encode a topology as a zero pattern and
ask whether the corresponding space of system matrices contains a Hurwitz
matrix. Belabbas addressed the difficulty above by bracketing stability between a
necessary and a sufficient graphical condition. 
For symmetric topologies,
\citep{kirkoryan2014decentralized,belabbas2022stable,belabbas2025structural}
showed that the two conditions coincide with stability, so that stability is
polynomial-time decidable. 
For general graphs, understanding how tight the bracket is, and the
complexity of deciding stability, remains an open problem; such is the starting point of the present work. 

There are two features that make the problem worth a meticulous study. First, the gap between the
two conditions is not a technicality of small examples: we show it is non-empty
on both sides and survives strong connectivity, so neither condition can be
promoted to a characterization by a local fix. The reduction to strongly
connected components, which underlies this part of the analysis, connects to
graph-augmentation techniques for enforcing strong connectivity in multi-agent
networks~\citep{ramos2025distributed}. Second, the two conditions turn
out to have opposite computational characters, one easy and one hard, which
sharpens the sense in which stability sits genuinely between them. 

\paragraph*{Main contributions.}
We show that the gap between the necessary and the sufficient condition of~\citep{belabbas2013sparse} is non-empty from both sides, and that it persists
even for strongly connected graphs. In particular, we exhibit a packable graph
that is not stable (of order four) and a stable graph that is not peelable (of
order three), refuting both converses of Belabbas' inclusions.

We then refine the necessary condition in two ways: first, by asking each
strongly connected component to be packable; second, by inspecting all entries
of the Routh-Hurwitz table~\citep{routh1877treatise,hurwitz1895bedingungen}
rather than only the coefficients of the characteristic polynomial. The second
refinement is natural: Belabbas' packability condition reads only whether the
coefficients of the characteristic polynomial can be made nonzero, whereas
stability requires the stronger Routh-Hurwitz inequalities among those
coefficients to hold simultaneously. Passing from the coefficients to the full
table captures algebraic relations, such as a forced identity among Hurwitz
determinants, that the cycle structure alone cannot see. We show that
both refinements are strict, exhibiting separating digraphs, and obtain the
chain of strict implications
\begin{center}
\begin{tikzpicture}[
  box/.style={draw, rounded corners=2pt, inner sep=3pt, font=\small, align=center},
  imp/.style={-{Stealth[length=2mm]}, thick}]
  \node[box] (a) {peelable};
  \node[box, right=6mm of a] (b) {stable};
  \node[box, right=6mm of b] (c) {SCC-Routh-Hurwitz};
  \node[box, right=6mm of c] (d) {SCC-packable};
  \node[box, right=6mm of d] (e) {SCC-rooted and packable};
  \draw[imp] (a) -- (b);
  \draw[imp] (b) -- (c);
  \draw[imp] (c) -- (d);
  \draw[imp] (d) -- (e);
\end{tikzpicture}
\end{center}

We next settle the complexity of the associated decision problems. We remark that packability and peelability have immediate NP certificates, and to the best of our knowledge, no other complexity bound
has previously been identified. 
Via a reduction to polynomial identity testing, each necessary condition (Belabbas'
packability and its refinements) lies in the class $\mathrm{RP}$ of
randomized polynomial time with one-sided error, hence is not $\mathrm{NP}$-hard
unless $\mathrm{NP}=\mathrm{RP}$. 
By contrast, a reduction
from satisfiability \citep{cook2023complexity} shows that the sufficient
condition (peelability) is $\mathrm{NP}$-complete; structural problems of this
kind can indeed be hard, as we have recently shown for the robust minimal
structural controllability problem \citep{pocas2026nphardness}. 
We also give the
corresponding algorithms: probabilistic tests for the necessary conditions and a
backtracking, worst-case exponential, procedure for peelability.

Finally, we implement these algorithms and report numerical experiments across
graph sizes (Section~\ref{sec:experiments}). For $n=5$, exhaustive enumeration over the $33\,554\,432$ digraphs shows that $21\,793\,941$ satisfy the necessary condition, $21\,787\,581$ are stable, and $20\,445\,031$ are peelable. 
For larger values of $n$, we estimate the probability that a random
digraph with $m$ edges satisfies the necessary and the sufficient
conditions. 
We report that the gap between the necessary and the sufficient
condition persists, and that the running times of our implementations grow
polynomially for the former and exponentially for the latter. 
The analytical and numerical results show that the tested combinatorial conditions, based on cycles and connectivity structures, do not provide a complete characterization of structural stability, an algebraic property. 
Whether a tractable graphical criterion exists remains an open problem.


\paragraph*{Organization.}
Section~\ref{sec:prelim} fixes notation and recalls the graphical conditions of
\cite{belabbas2013algorithms,belabbas2013sparse}. Section~\ref{sec:scc} reduces the study to strongly
connected patterns and introduces the SCC-packable condition.
Section~\ref{sec:routh-hurwitz} introduces the stronger Routh-Hurwitz necessary
condition, and Section~\ref{sec:strict} establishes the chain of strict
inclusions between the conditions through separating digraphs.
Section~\ref{sec:complexity} places the necessary conditions in $\mathrm{RP}$,
and Section~\ref{sec:sufficient-hard} proves the sufficient condition
$\mathrm{NP}$-complete. Section~\ref{sec:experiments} reports the numerical
experiments.

\section{Preliminaries and notation}\label{sec:prelim}

Let $n>0$ and $\alpha\subseteq\{1,\dots,n\}^2$. The associated \emph{sparse
matrix space} (SMS) is
\begin{equation}
\Sigma_\alpha=\{A\in\R^{n\times n}: a_{uv}=0 \text{ for } (u,v)\notin\alpha\};
\end{equation}
the parameters indexed by $\alpha$ are free. To $\Sigma_\alpha$ we attach a
digraph $G=(V,E)$ with $V=\{1,\dots,n\}$, $E\subseteq V\times V$, and an edge
$(u,v)\in E$ exactly when $a_{uv}$ is free; a self-loop $(u,u)$ corresponds to a
free diagonal parameter and its endpoint is called \emph{rooted}. We write $\SG$
for the SMS of $G$ and call $G$ (and $\SG$) \emph{stable} when $\SG$ contains a
Hurwitz matrix.

A \emph{subgraph} $H=(V',E')$ of a digraph $G=(V,E)$ is a digraph with
$V'\subseteq V$ and $E'\subseteq E\cap(V'\times V')$. A \emph{cycle} is a
sequence of edges $(x_0,x_1),(x_1,x_2),\dots,(x_{k-1},x_0)$ with the $x_i$
distinct; its \emph{size} is its number of edges, and a self-loop is a cycle of
size one. A \emph{cycle packing} (a.k.a. \textbf{Hamiltonian decomposition}~\citep{belabbas2013sparse}) is a set of vertex-disjoint cycles; its size is
the sum of the sizes of its cycles, equal to the number of vertices it covers. A
\emph{path} of size $k$ is a sequence of distinct vertices $x_0,x_1,\dots,x_k$ such that $(x_0,x_1),\dots,(x_{k-1},x_k)$ are edges in the digraph;
a single vertex is a trivial path of size $k=0$. A digraph is \emph{strongly
connected} if there is a path between every ordered pair of vertices, and a
\emph{strongly connected component} (SCC) is a maximal strongly connected
subgraph.

Following \cite{belabbas2013sparse}, we consider the following graph properties
related to structural stability; the terminology \emph{packable} and
\emph{peelable} is ours.
\begin{itemize}
\item $G$ is \emph{SCC-rooted} if every strongly connected component of $G$ has
a rooted vertex.
\item $G$ is \emph{packable} if $G$ has a cycle packing of size $k$ for every
$k\in\{1,\dots,n\}$.
\item $G$ is \emph{peelable} if it contains a nested sequence
$G_1\subset G_2\subset\cdots\subset G_n=G$ in which each $G_k$ has a cycle
packing of size $|V(G_k)|=k$.
\end{itemize}

The following theorem records the two graphical conditions of
\cite{belabbas2013sparse} that frame the paper.
\begin{theorem}[\cite{belabbas2013sparse}]\label{thm:belabbas}
Every stable digraph is SCC-rooted and packable. Furthermore, every peelable
digraph is stable.\hfill$\diamond$
\end{theorem}

Theorem~\ref{thm:belabbas} leaves open the gap between the necessary and the
sufficient condition and the complexity of the associated decision problems,
which are the subject of this paper.

\section{Reduction to strongly connected patterns}\label{sec:scc}

As noted in \citet{belabbas2022stable}, we may split the stability question of a digraph into its strongly connected components.

\begin{lemma}\label{lem:scc}
Let $G$ be a digraph with strongly connected components $ S_1,\ldots, S_\ell$. Then $G$ is stable \textit{iff} each of $S_1,\ldots,S_\ell$ is stable.~\hfill$\diamond$
\end{lemma}


\begin{myproof}
Sort $S_1,\ldots S_\ell$ by their topological order. After
the corresponding symmetric permutation, every member $A$ of the matrix space $\SG$ is block upper triangular, since an
edge between distinct components runs forward in the order. Moreover, the diagonal blocks $A_1,\ldots, A_\ell$ belong in the matrix spaces of $S_1,\ldots,S_\ell$ respectively, and the spectrum of $A$ is
the union of the spectra of the diagonal blocks, so $G$ is stable if and
only if there exist a Hurwitz matrix $A\in\SG$, if and only if there exist Hurwitz matrices $A_i\in\Sigma_{S_i}$, if and only if each $S_i$ is stable. 
\end{myproof}

The practical relevance of Lemma~\ref{lem:scc} is that partitioning a graph $G$ into its SCCs can be done in polynomial time by Tarjan's algorithm~\citep{tarjan1972depth}, enabling a divide-and-conquer approach to study the original problem. A natural question is whether Belabbas' conditions also decompose into the digraph SCCs. For peelability, the answer is yes.

\begin{lemma}\label{lem:peel-scc}
Let $G$ be a digraph with strongly connected components $ S_1,\ldots, S_\ell$. Then $G$ is peelable \textit{iff} each of $S_1,\ldots,S_\ell$ is peelable.~\hfill$\diamond$
\end{lemma}

\begin{myproof}
Note that every cycle lies inside a single strongly connected component. Hence a cycle packing of $G$ can be partitioned into cycle packings on each component $S_i$, and conversely, cycle packings for each component assemble into one of $G$. 

Suppose that $G$ is peelable with the nested sequence $G_1\subset G_2\subset\cdots\subset G_n=G$, such that each $G_k$ has a cycle packing of size $|V(G_k)|=k$. For each SCC $S_i$, let $S_{i,k} = S_i\cap G_k$. By the above discussion, each $S_{i,k}$ has a cycle packing of size $|V(S_{i,k})|$; when increasing $k$, the size $|V(S_{i,k})|$ increases by one or remains the same depending on whether the vertex added is in $S_i$. By taking an appropriate subsequence of $S_{i,k}$, we conclude that $S_i$ is peelable. 

Conversely, suppose that each $S_i$ is peelable. By listing the strongly connected components in any order and following the nested sequence within each component, we obtain a nested sequence of subgraphs of $G$ with cycle packings of each size. Therefore, $G$ is peelable.
\end{myproof}

For Belabbas' packability condition, the following Lemma shows that packability of each component implies packability of the original graph. In Section~\ref{sec:strict}, we provide a counter-example for the converse. 

\begin{lemma}\label{lem:scc-packable}
Let $G$ be a digraph with strongly connected components $ S_1,\ldots, S_\ell$. If each $S_i$ is packable, then so is~$G$.~\hfill$\diamond$
\end{lemma}

\begin{myproof}
Suppose that each $S_i$ is packable. For each $k\in\{1,\ldots,n\}$, take any choice $k=k_1+\ldots+k_\ell$ such that $0\leq k_i\leq |V(S_i)|$. The union of $\ell$ cycle packings, one for each $S_i$ of size $k_i$, produces a cycle packing for $G$, of size $k$. Therefore, $G$ is packable. 
\end{myproof}
    
The above discussion suggests a stronger necessary condition for stability.

\begin{itemize}
\item $G$ is \emph{SCC-packable} if every strongly connected component of $G$ is packable.
\end{itemize}

\begin{lemma}\label{lem:scc-packable-necessary}
Every stable digraph is SCC-packable. Furthermore, every SCC-packable digraph is SCC-rooted and packable.~\hfill$\diamond$
\end{lemma}

\begin{myproof}
Suppose that $G$, with strongly connected components $S_1,\ldots,S_\ell$, is stable. Then, by Lemma~\ref{lem:scc}, each $S_i$ is stable, and by Theorem~\ref{thm:belabbas}, each $S_i$ is packable, so that $G$ is SCC-packable. 

Suppose now that $G$ is SCC-packable, that is, each $S_i$ is packable. By taking $k=1$, each $S_i$ has a cycle packing of size 1, i.e., each $S_i$ has a rooted vertex; thus, $G$ is SCC-rooted. Moreover, $G$ is packable by Lemma~\ref{lem:scc-packable}.
\end{myproof}

\section{A stronger necessary condition via Routh-Hurwitz tables}\label{sec:routh-hurwitz}

As already pointed out by Belabbas, the existence of cycle packings on a digraph $G$ is tightly related with the characteristic polynomial coefficients of matrices in $\SG$. 
For $A\in\R^{n\times n}$, we denote its characteristic polynomial~by
\begin{equation}\label{eq:char-poly}
\chi_A(s)=\det(sI-A)=s^{n}+p_1 s^{n-1}+p_2 s^{n-2}+\cdots+p_n .
\end{equation}

The following characterization is a graph-theoretic formulation of the determinant expansion in~\cite{kaplan1962operational}; see also~\cite{belabbas2013sparse}.

\begin{lemma}[\cite{kaplan1962operational}]\label{lem:pk-packing}
Let $G$ be a digraph with $n$ vertices and $k\in\{1,\ldots,n\}$. Then $G$ has a cycle packing of size $k$ if and only if there is a matrix $A\in\Sigma_G$ whose characteristic polynomial satisfies $p_k\neq 0$.~\hfill$\diamond$
\end{lemma}

\begin{myproof}
By the cycle-family expansion of the determinant, $p_k=\sum_{L\in\mathcal{L}_k}(-1)^{c(L)}w(L)$,
where $\mathcal{L}_k$ is the family of cycle packings of size $k$, and for each cycle packing $L\in\mathcal{L}_k$,
$c(L)$ is the number of cycles of $L$ and $w(L)$ the product of the entries along them. 
Distinct packings use distinct edge sets and hence yield distinct monomials of the same degree $k$, which cannot cancel. 
So $p_k$ is identically zero (over every $A\in\SG$) if and only if $\mathcal{L}_k=\emptyset$, that is, if and only if $G$ has no cycle packing of size $k$.
\end{myproof}

We refine the packability condition by inspecting the full Routh-Hurwitz table~\citep{routh1877treatise,hurwitz1895bedingungen}, which includes additional relations among the polynomial coefficients going beyond cycle structure analysis.

\begin{definition}\label{def:rh-table}
Let $A\in\R^{n\times n}$ with characteristic polynomial as in \eqref{eq:char-poly}. The \emph{Routh-Hurwitz table} of $A$ is the array 
\[
\begin{array}{>{$}p{18mm}<{$}>{$}p{18mm}<{$}>{$}p{18mm}<{$}l}
R_{0,0}= 1 & R_{0,1}= p_2 & R_{0,2}= p_4 & \cdots\\
R_{1,0}= p_1 & R_{1,1}= p_3 & R_{1,2}= p_5 & \cdots\\
R_{2,0} & R_{2,1} & R_{2,2} & \cdots\\
R_{3,0} & R_{3,1} & R_{3,2} & \cdots\\
\vdots & \vdots & \vdots & \ddots
\end{array}
\]
whose first two rows are
the alternate polynomial coefficients
and whose subsequent rows are formed by the standard recursion: if two consecutive rows are $(R_{i-2,0},R_{i-2,1},\dots)$ and $(R_{i-1,0},R_{i-1,1},\dots)$, with $R_{i-1,0}\neq 0$, then the next row is
$(R_{i,0},R_{i,1},\dots)$ with
\begin{equation}
R_{i,j}=R_{i-2,j+1}-\frac{R_{i-2,0}}{R_{i-1,0}}R_{i-1,j+1}.
\label{eq:rh-recursion}
\end{equation}
The table is \emph{well-defined} if none of the denominators in~\eqref{eq:rh-recursion} vanish.~\hfill$\circ$
\end{definition}

We denote by $R_{i,j}$ the $(i,j)$ entry of the table, for $i+2j\leq n$; it is not hard to see that there are exactly $\lfloor(n+2)^2/4\rfloor = \mathcal{O}(n^2)$ such entries. 
It is known that, whenever $i+2j > n$, the entries $R_{i,j}$ are identically zero, and we call such entries \emph{structural zeros}~\citep{routh1877treatise} and the remaining entries \emph{structural positions}. 
 The entries of the \emph{first column} of the table ($R_{0,0}=1,R_{1,0}=p_1,R_{2,0},\ldots,R_{n,0})$ are, up to normalization, the leading principal minors of the Hurwitz matrix
 $H=[\,p_{2c-r}\,]_{r,c=1}^{n}$ (with $p_0=1$ and $p_k=0$ for $k\notin\{0,\dots,n\}$).
Writing $\Delta_0=1$ and $\Delta_k=\det H_{[1..k],[1..k]}$ for the $k$-th Hurwitz determinant, the first-column
entries are the successive quotients 
$R_{0,0}=\Delta_0, R_{k,0}=\Delta_k/\Delta_{k-1}$ $(1\le k\le n)$, whenever the denominators do not vanish. 

\begin{lemma}[Routh's criterion~\citep{routh1877treatise}]\label{lem:routh}
A matrix $A\in\R^{n\times n}$ is Hurwitz if and only if every entry in the first column of its Routh-Hurwitz table is well-defined and positive, equivalently, $\Delta_k>0$ for all $k\in\{1,\ldots,n\}$.~\hfill$\diamond$ 
\end{lemma}

Recall that, for a digraph $G$, we may think of a generic matrix $A\in\SG$ with free parameters $a_{uv}$ indexed by the edges of $G$, $(u,v)\in E$. By the cycle-family expansion of the determinant, each of the coefficients $p_k$ can be understood as a homogeneous polynomial of degree $k$ in the parameters $a_{uv}$. Furthermore, the Routh-Hurwitz table entries $R_{i,j}$ can be seen as rational functions (quotients of polynomials) in the free entries. 
A meaningful strengthening of packability is, thus, to ask whether each Routh-Hurwitz table entry (other than a structural zero) is identically zero over $\SG$ or non-zero for some choices of parameters. 

\begin{definition}\label{def:scc-rh}
A digraph $G$ is \emph{Routh-Hurwitz} if there is a matrix $A\in\SG$ whose Routh-Hurwitz table is well-defined and satisfies $R_{i,j}\neq 0$ for each structural position ($i+2j\leq n)$. 
A digraph $G$ is \emph{SCC-Routh-Hurwitz} (SCC-RH) if every strongly connected component of $G$ is \emph{Routh-Hurwitz}.~\hfill$\circ$
\end{definition}

We shall observe in Section~\ref{sec:complexity} that the above graph property can be checked efficiently, by reducing the problem to polynomial identity testing, decidable in randomized
polynomial time. 
For now, we formally compare our property with stability and packability.

\begin{theorem}\label{thm:stable-implies-rh-implies-packable}
Every stable digraph is SCC-Routh-Hurwitz. Furthermore, every SCC-Routh-Hurwitz digraph is SCC-packable.~\hfill$\diamond$
\end{theorem}



\begin{myproof}
Let $G$ be stable, and let $S$ be any strongly connected component of $G$, with
$n=|V(S)|$. By Lemma~\ref{lem:scc}, $S$ is stable, so some $A^\ast\in\Sigma_S$ is
Hurwitz. By Routh's criterion (Lemma~\ref{lem:routh}), the first-column entries of
its Routh-Hurwitz table are positive, $R_{i,0}>0$ for $0\leq i\leq n$. Rewriting the recursion \eqref{eq:rh-recursion} as
\begin{equation}\label{eq:rh-recursion-backwards}
R_{i,j+1} = R_{i+2,j} + \frac{R_{i,0}}{R_{i+1,0}}\,R_{i+1,j+1},
\end{equation}
we argue by induction on $i$ decreasing from $n$ to $0$. The two bottom rows have positive first-column entries by assumption; in the inductive step, looking at \eqref{eq:rh-recursion-backwards}, the term
$R_{i+2,j}$ is positive, the term $\frac{R_{i,0}}{R_{i+1,0}}$ is positive, and the entry $R_{i+1,j+1}$ is nonnegative (it is either positive or a structural zero). Hence, $R_{i,j}>0$ for all $i+2j\leq n$.
In particular, no entry at a structural position vanishes, so $S$ is Routh-Hurwitz. As $S$ was
arbitrary, $G$ is SCC-RH.

Next, suppose that $G$ is SCC-RH, and again let $S$ be a component
with $n=|V(S)|$. Since the coefficients $p_1,\dots,p_n$ of the characteristic polynomial
over $\Sigma_S$ are the entries in the top two rows of the Routh-Hurwitz table, and $S$ is Routh-Hurwitz, we conclude that no $p_k$ is identically zero. 
By Lemma~\ref{lem:pk-packing}, $S$ has a
cycle packing of every size $k\in\{1,\dots,n\}$, that is, $S$ is packable. This holds for any component $S$, and thus $G$ is SCC-packable.~
\end{myproof}



\section{A chain of strict inclusions}\label{sec:strict}

By Theorem~\ref{thm:stable-implies-rh-implies-packable}, the SCC-RH condition sits
between stability and SCC-packability. 
In this section, we provide examples that separate these properties, as well as peelability (the sufficient condition for stability) and SCC-rooted packability (Belabbas' original necessary condition).

\begin{proposition}\label{prop:packable-not-scc-packable}
There exists an SCC-rooted and packable digraph that is not SCC-packable.~\hfill$\diamond$
\end{proposition}

\begin{figure}[!ht]
\centering
\begin{tikzpicture}[scale=.6, transform shape,node distance=1.5cm]
\begin{scope}[every node/.style={circle,thick,draw},square/.style={regular polygon,regular polygon sides=4}]
\node[label={[label distance=-1mm] below:{\Large $1$}}] (1) at (0.,1.05868) {};
\node[label={[label distance=-1mm] below:{\Large $2$}}] (2) at (2.78621,1.05864) {};
\node[label={[label distance=-1mm] right:{\Large $3$}}] (3) at (5.07893,2.11728) {};
\node[label={[label distance=-1mm] right:{\Large $4$}}] (4) at (5.07883,0.) {};
\end{scope}
\begin{scope}[>={Stealth[black]},
              every edge/.style={draw=black, thick}]
\path [->] (1) edge[loop above] node {} (1);
\path [->] (1) edge node {} (2);
\path [->] (2) edge[loop above] node {} (2);
\path [->] (2) edge node {} (3);
\path [->] (3) edge node {} (4);
\path [->] (4) edge node {} (2);
\end{scope}
\end{tikzpicture}
\caption{A packable digraph with a strongly connected component that is not packable.}
\label{fig:inpackable}
\end{figure}
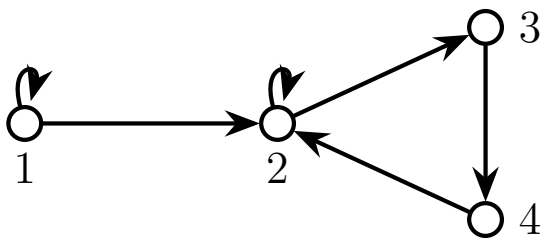

\begin{myproof}
Consider the graph in Fig.~\ref{fig:inpackable}. The graph admits the cycle packings on vertex sets $\{1\}$ (self-loop), $\{1,2\}$ (two self-loops), $\{2,3,4\}$ (3-cycle) and $\{1,2,3,4\}$ (self-loop plus 3-cycle). However, the SCC on vertices $\{2,3,4\}$ is not packable, since it does not contain a cycle packing of size 2.
\end{myproof}

\begin{proposition}\label{prop:packable-not-rh}
There exists a strongly connected digraph that is packable but not Routh-Hurwitz.~\hfill$\diamond$
\end{proposition}

\begin{figure}[!ht]
\centering
\begin{tikzpicture}[scale=.6, transform shape,node distance=1.5cm]
\begin{scope}[every node/.style={circle,thick,draw},square/.style={regular polygon,regular polygon sides=4}]
\node[label={[label distance=-1mm] above:{\Large $1$}}] (1) at (0.00612201,-0.00612201) {};
\node[label={[label distance=-1mm] left:{\Large $4$}}] (2) at (1.77609,-1.76865) {};
\node[label={[label distance=-1mm] right:{\Large $3$}}] (3) at (3.54021,0.00741088) {};
\node[label={[label distance=0mm] left:{\Large $2$}}] (4) at (1.76731,1.76731) {};
\end{scope}
\begin{scope}[>={Stealth[black]},
              every edge/.style={draw=black, thick}]
\path [->] (1) edge[loop left] node {} (1);
\path [->] (1) edge node {} (4);
\path [->] (2) edge[bend right=15] node {} (3);
\path [->] (3) edge[bend right=15] node {} (2);
\path [->] (3) edge[bend right=15] node {} (4);
\path [->] (4) edge[bend right=15] node {} (3);
\path [->] (2) edge node {} (1);
\end{scope}
\end{tikzpicture}

\caption{Example of a strongly connected digraph that is packable but not Routh-Hurwitz.}
\label{fig:g4}
\end{figure}
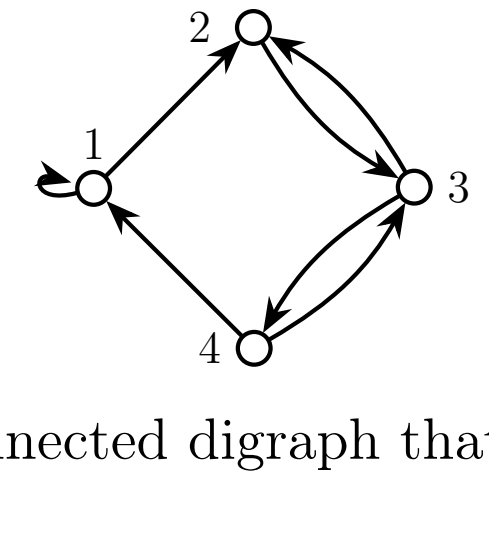

\begin{myproof}
Consider the strongly connected digraph $G$ on 4 vertices depicted in
Fig.~\ref{fig:g4} (Example 3.5 in~\cite{cavers2021polynomial} identifies this digraph as unstable). The graph admits the cycle packings on vertex sets $\{1\}$ (self-loop), $\{2,3\}$ (2-cycle), $\{1,2,3\}$ (self-loop plus 2-cycle) and $\{1,2,3,4\}$ (4-cycle), so that $G$ is packable. Its characteristic polynomial has coefficients $p_1=-a_{11}$,
$p_2=-a_{23}a_{32}-a_{34}a_{43}$, $p_3=a_{11}a_{23}a_{32}+a_{11}a_{34}a_{43}=p_1p_2$, $p_4=-a_{12}a_{23}a_{34}a_{41}$. In particular, in the Routh-Hurwitz table, we get
$
R_{2,0}=p_2-\frac{1}{p_1}p_3=p_2-p_2\equiv 0 .
$ 
This table entry is either undefined (whenever $p_1$ vanishes) or zero over $\Sigma_{G}$, so $G$ is not
Routh-Hurwitz. 
\end{myproof}

\begin{proposition}\label{prop:rh-not-sufficient}
There exists a strongly connected digraph that is Routh-Hurwitz but not
stable.~\hfill$\diamond$
\end{proposition}

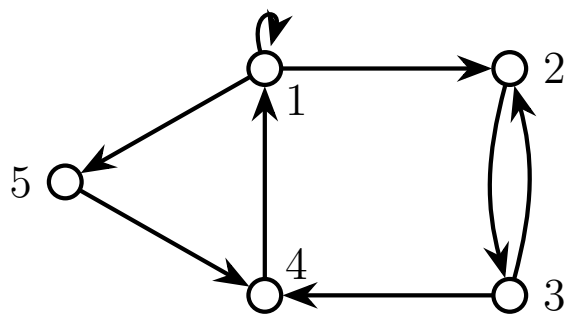
\begin{figure}[!ht]
    \centering
    \begin{tikzpicture}[scale=.6, transform shape,node distance=1.5cm]
\begin{scope}[every node/.style={circle,thick,draw},square/.style={regular polygon,regular polygon sides=4}]
\node[label={[label distance=-1mm] right:{\Large $3$}}] (3) at (0,0) {};
\node[label={[label distance=-1mm] right:{\Large $2$}}] (2) at (0,2.5) {};
\node[label={[label distance=-1mm] below right:{\Large $1$}}] (1) at (-2.7,2.5) {};
\node[label={[label distance=-1mm] above right:{\Large $4$}}] (4) at (-2.7,0) {};
\node[label={[label distance=-1mm] left:{\Large $5$}}] (5) at (-4.9,1.25) {};
\end{scope}
\begin{scope}[>={Stealth[black]},
              every edge/.style={draw=black, thick}]
\path [->] (2) edge[bend right=15] node {} (3);
\path [->] (3) edge[bend right=15] node {} (2);
\path [->] (1) edge[loop above] node {} (1);
\path [->] (4) edge node {} (1);
\path [->] (1) edge node {} (5);
\path [->] (5) edge node {} (4);
\path [->] (1) edge node {} (2);
\path [->] (3) edge node {} (4);

\end{scope}
\end{tikzpicture}

    \caption{Example of a digraph that is SCC-RH but that is not stable.}
    \label{fig:scc_r_h_unstable}
\end{figure}

\begin{myproof}
Consider the strongly connected digraph $G$ on 5 vertices depicted in Fig.~\ref{fig:scc_r_h_unstable}. Take, for example, the matrix $A\in\SG$ with $a_{11}=-2,a_{23}=1$, and the remaining free entries equal to $-1$. In this way, the Routh-Hurwitz table of $\chi_A$ is
\[
\def\linecompress{-2mm}
\begin{array}{>{$}p{15mm}<{$}>{$}p{15mm}<{$}>{$}p{15mm}<{$}}
 1 & 1 & 1 \\
 2 & 3 & 1 \\
 -0.5 & 0.5 & 0 \\
 5 & 1 & 0 \\
 0.6 & 0 & 0 \\
 1 & 0 & 0 \\
\end{array}
\]
so that no structural position vanishes. We conclude that $G$ is Routh-Hurwitz. 
Now suppose that $G$ was stable, with some Hurwitz matrix $A$. 
Write $c=a_{11}$, $u=a_{23}a_{32}$
for the $2$-cycle product on $\{2,3\}$, and $v=a_{15}a_{54}a_{41}$ for the
$3$-cycle product on $\{1,4,5\}$. The conditions $p_1>0$, $p_2>0$ and $p_5>0$ can be rewritten as $c<0$, $u<0$ and $uv>0$ respectively, which also imply $v<0$. However, these conditions force $R_{2,0}=-v/c<0$, a contradiction.
Therefore $\Sigma_{G}$ contains no Hurwitz matrix and $G$ is unstable.
\end{myproof}

\begin{proposition}\label{prop:notnec}
There exists a strongly connected digraph that is stable but not peelable.~\hfill$\diamond$
\end{proposition}

\begin{figure}[!ht]
\centering
\begin{tikzpicture}[scale=.6, transform shape,node distance=1.5cm]
\begin{scope}[every node/.style={circle,thick,draw},square/.style={regular polygon,regular polygon sides=4}]
  \node[label={[label distance=-1mm] above:{\Large $1$}}] (v1) at (0,0) {};
  \node[label={[label distance=-1mm] right:{\Large $2$}}] (v2) at (1.4,1) {};
  \node[label={[label distance=-1mm] right:{\Large $3$}}] (v3) at (1.4,-1) {};
  \end{scope}
  \begin{scope}[>={Stealth[black]},
              every edge/.style={draw=black, thick}]
  \draw[->] (v1) edge (v2);
  \draw[->] (v2) edge[bend left=18] (v3);
  \draw[->] (v3) edge[bend left=18] (v2);
  \draw[->] (v3) edge (v1);
  \draw[->] (v1) edge[loop left] (v1);
  \end{scope}
\end{tikzpicture}
\caption{A strongly connected digraph $G$ which is stable but not peelable.}
\label{fig:g3}
\end{figure}
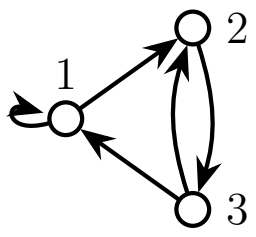

\begin{myproof}
Consider the strongly connected digraph $G$ on 3 vertices depicted in Fig.~\ref{fig:g3} (Fig.~4 in~\citep{cavers2021polynomial} identifies this digraph as stable). Take, for example, the matrix $A\in\SG$ with $a_{11}=a_{32}=-1,a_{12}=a_{23}=1,a_{31}=1/2$. The characteristic polynomial $\chi_A(s)=s^3+s^2+s+1/2$ has roots $-0.65$
and $-0.18\pm 0.86\,\mathrm{i}$, hence $A$ is Hurwitz and $G$ is stable.

On the other hand, the only cycle packing of size $1$ is the self-loop in $\{1\}$ and the only cycle packing of size $2$ is the 2-cycle in $\{2,3\}$. 
Hence, there are no subgraphs $G_1\subset G_2$ of sizes $1$ and $2$ having full cycle packings, implying that $G$ is not peelable.
\end{myproof}

To summarize, we have established the following chain of strict inclusions between graph properties related to stability, with two refinement steps between stability and Belabbas' necessary condition of packability:
\[
\resizebox{\textwidth}{!}{$
\text{peelable}\ 
\overset{\text{Thm.~\ref{thm:belabbas}}}{\underset{\text{Prop.~\ref{prop:notnec}}}{\Longrightarrow}} 
\text{stable}\ 
\overset{\text{Thm.~\ref{thm:stable-implies-rh-implies-packable}}}{\underset{\text{Prop.~\ref{prop:rh-not-sufficient}}}{\Longrightarrow}}\ 
\text{SCC-Routh-Hurwitz}\ 
\overset{\text{Thm.~\ref{thm:stable-implies-rh-implies-packable}}}{\underset{\text{Prop.~\ref{prop:packable-not-rh}}}{\Longrightarrow}}\ 
\text{SCC-packable}\ 
\overset{\text{Lem.~\ref{lem:scc-packable-necessary}}}{\underset{\text{Prop.~\ref{prop:packable-not-scc-packable}}}{\Longrightarrow}}\ 
\text{SCC-rooted and packable}.$}
\]

\begin{remark}
Notice that adding an isolated vertex with a self-loop to a graph preserves all properties of interest in both directions. 
Thus, the examples above also yield separations at every larger order.
In other words, the gaps between peelability, stability, and each of the necessary conditions persist at every number of vertices beyond five.~\hfill$\nabla$
\end{remark}

\section{The necessary conditions are in RP}\label{sec:complexity}

It is natural to ask how hard it is to test the necessary and sufficient conditions of \cite{belabbas2013sparse}, which are combinatorial in nature.
In this section, we show that the necessary condition (existence of cycle packings for every size $k$), as well as both our refined necessary conditions, are decidable by probabilistic polynomial-time algorithms (with one-sided errors), placing them in the complexity class $\mathrm{RP}\subseteq\mathrm{NP}$~\cite[Chapter 7]{arora2009computational}.
Consequently, none of these conditions is $\mathrm{NP}$-hard, unless $\mathrm{NP}=\mathrm{RP}$.

Let us provide some intuition to see why a probabilistic algorithm exists. Suppose that $G$ has a cycle packing of size $k$. By Lemma~\ref{lem:pk-packing}, the coefficient $p_k$ of the characteristic polynomial is not identically zero.
Crucially, we can view $p_k$ as a polynomial expression in the matrix entries $a_{ij}$. Its zero set is therefore an algebraic surface of Lebesgue measure zero.
Thus, almost every assignment of the entries is a witness for $p_k\neq 0$, assuming a cycle packing of size $k$ exists.

To turn this intuition into a formal algorithm, and at the same time avoid potential floating-point inaccuracies, we evaluate the coefficients at random integer points and use the Schwartz-Zippel Lemma~\cite[Lemma A.36]{arora2009computational} to bound the probability of a false negative. For a positive integer $K$, we write $A\sim\Sigma_G(K)$ to mean that $A$ is obtained by setting $a_{uv}=0$ whenever $(u,v)\notin E$, and sampling each remaining entry independently and uniformly from $\{1,\ldots,K\}$. As before, we denote by $p_k(A)$ the value of the $k$-th coefficient of the characteristic polynomial $\chi_A$, as a function of (the entries of) $A$. %

\begin{proposition}\label{prop:computation-pk}
Given an $n\times n$ integer matrix $A$ with entries in $\{0,\ldots,K\}$, the values $p_k(A)$ for $k=1,\ldots,n$ can be computed exactly in time polynomial in $n$ and $\log(K)$.~\hfill$\diamond$
\end{proposition}

\begin{myproof}
The Faddeev-LeVerrier algorithm~\citep{hou1998classroom} computes the coefficients
of $\chi_A$ by the recursion
\begin{equation}\label{eq:faddeev}
A_1=A,\qquad p_k=-\frac{1}{k}\tr(A_k),\qquad A_{k+1}=A(A_k+p_kI),
\end{equation}
in which every division by $k$ is exact over the integers. We now bound the absolute values of the entries of $A_k$, assuming $\|A\|_\infty\coloneqq\max|a_{ij}|\leq K$:
\begin{equation*}
\begin{array}{rcl}
\|A_{k+1}\|_\infty&=&\|A(A_k+p_kI)\|_\infty \leq n\|A\|_\infty \|A_k+p_kI\|_\infty\\
& \leq & nK(\|A_k\|_\infty+|p_k|)
\leq nK(\|A_k\|_\infty+n\|A_k\|_\infty)\\
& \leq & n(1+n)K\|A_k\|_\infty
\leq 2n^2K\|A_k\|_\infty.
\end{array}
\end{equation*}
Therefore $\|A_k\|_\infty \leq (2n^2K)^k$. Taking logarithms, every integer produced has $\mathcal O(n\log(nK))$ bits. Each iteration computes one matrix product and one trace, so the total running time is polynomial in $n$ and $\log(K)$.
\end{myproof}

\begin{lemma}\label{lem:pk-probability-bound}
Let $G$ be a digraph on $n$ vertices and $k\in\{1,\ldots,n\}$.
\begin{enumerate}
\item If $p_k$ is identically zero, then $\mathbb{P}_{A\sim\Sigma_G(K)}[p_k(A)=0]=1$;
\item If $p_k$ is not identically zero, then $\mathbb{P}_{A\sim\Sigma_G(K)}[p_k(A)=0]\leq k/K$.~\hfill$\diamond$
\end{enumerate}
\end{lemma}

\begin{myproof}
The first item comes from the definition of identically zero polynomial. For the second, $p_k$ is a nonzero polynomial of degree $k$, so the bound follows from the Schwartz-Zippel Lemma~\cite[Lemma A.36]{arora2009computational}. 
\end{myproof}

\begin{theorem}\label{thm:cycle-packing-rp}
Packability of a given graph $G$ is decidable in $\mathrm{RP}$ (probabilistic polynomial time with one-sided errors). Formally, there is a probabilistic polynomial-time algorithm that receives a digraph $G$ and: returns \true\ with probability at least $1/2$, if $G$ is packable; returns \false\ with probability $1$, if $G$ is not packable.~\hfill$\diamond$
\end{theorem}

\begin{myproof}
Consider Algorithm~\ref{alg:packability}, which receives a digraph $G$ as well as additional parameters $K$ (sampling range) and $T$ (number of trials).

\begin{algorithm}[!ht]
\caption{Randomized Packability Test}
{\small
\begin{algorithmic}[1]
\Require digraph $G$ with $n$ vertices, integer $K\geq 2n$, number of trials $T$
\For{$k=1$ to $n$}
$r_k \gets \false$
\Comment{{\color{gray}\scriptsize initialization}}
\EndFor
\For{$t=1$ to $T$}
\Comment{{\color{gray}\scriptsize repeat $T$ times}}
\State $A\gets$ random matrix in $\Sigma_{G}(K)$ \label{line:start-trial}
\For{$k=1$ to $n$}
\If{$p_k(A)\neq 0$} $r_k \gets \true$\label{line:pk-zero}
\EndIf
\EndFor \label{line:end-trial}
\EndFor
\State\Return{whether every $r_k = \true$} \label{line:pk-check}
\end{algorithmic}
}
\label{alg:packability}
\end{algorithm}

Suppose that $G$ is not packable.
Then, there exists a $k\in\{1,\ldots,n\}$ such that $G$ has no cycle packings of size $k$.
Therefore, $p_k$ is identically zero, so that (by Lemma~\ref{lem:pk-packing}) $p_k(A)=0$ for every integer matrix $A\in\Sigma_G(K)$.
We conclude that variable $r_k$ remains \false\ throughout all trials, so that the algorithm correctly returns \false.

Next, suppose that $G$ is packable.
For each trial, and each $k\in\{1,\ldots,n\}$, the probability that $p_k(A)=0$ (Line~\ref{line:pk-zero}) is bounded by $k/K\leq n/2n=1/2$, thanks to Lemma~\ref{lem:pk-probability-bound}.
For each $k\in\{1,\ldots,n\}$, the probability that $r_k$ remains \false\ at Line~\ref{line:rk-check} equals the probability that $p_k(A)=0$ over all $T$ trials, which is bounded by $2^{-T}$ since trials are independent. By a union bound, the probability that at least one $r_k$ remains \false\ at Line~\ref{line:rk-check} is at most $n\cdot 2^{-T}$.
By taking $T\geq\lceil\log_2(2n)\rceil$, we ensure that the probability of a false negative is bounded by~$1/2$, as desired. Hence, only logarithmically many trials are required.

To conclude, notice that the total complexity of Algorithm~\ref{alg:packability} is at most $T$ times the complexity of a single trial (Lines~\ref{line:start-trial}--\ref{line:end-trial}), which itself runs in time polynomial in $n$, using Proposition~\ref{prop:computation-pk} and taking for example $K=2n$. 
\end{myproof}


\begin{corollary}\label{cor:scc-packable-rp}
Determining whether a graph $G$ is SCC-rooted and packable can be decided in $\mathrm{RP}$. Determining whether a graph $G$ is SCC-packable can be decided in~$\mathrm{RP}$.~\hfill$\diamond$
\end{corollary}

\begin{myproof}
To determine whether $G$ is SCC-rooted and packable, we begin by computing its strongly connected components, which can be done in polynomial time. 
If there is one component without a self-loop, we immediately reject $G$.
Otherwise, we run Algorithm~\ref{alg:packability} on $G$, with probability of a false negative at most $1/2$. 

To determine whether $G$ is SCC-packable, we again begin by computing its strongly connected components $S_1,\ldots,S_\ell$. 
We run $T'\geq \lceil\log_2(2\ell)\rceil$ independent trials of Algorithm~\ref{alg:packability} on each of the strongly connected components, which reduces the probability of error from $1/2$ to $2^{-T'}\leq 1/(2\ell)$. 
We accept $G$ if and only if at least one trial for each strongly connected component accepts. By a union bound, the final probability of error can be kept below~$1/2$. 
\end{myproof}


For the stronger necessary condition (SCC-RH), the table entries $R_{i,j}$ of Definition~\ref{def:rh-table} are no longer polynomials in the matrix entries $a_{ij}$ but rather rational functions. 
Their numerators, however, are polynomials: these are the bordered minors~\citep{Jury1974inners} of the Hurwitz matrix. 
Moreover, these minors have degree $\mathcal O(n^2)$, keeping the integer arithmetic feasible in practice.

\begin{definition}[\cite{Jury1974inners}]\label{def:bordered-minors}
Let $A\in\mathbb R^{n\times n}$ with characteristic polynomial as in~\eqref{eq:char-poly}. Let $H=\left[\,p_{2c-r}\,\right]_{r,c=1}^{n}$ (with $p_0=1$ and $p_k=0$ for $k\notin\{0,\ldots,n\}$) be the Hurwitz matrix of $A$. For a pair of nonnegative integers $(i,j)$ with $i+2j\leq n$, the \emph{bordered Hurwitz minor} is defined as:
\begin{equation}\label{eq:bordered-minor}
B_{0,j} = p_{2j}, \quad
B_{i,j} = \det H_{\{1,\ldots,i\},\{1,\ldots,i-1\}\cup\{i+j\}}\ (i\geq 1).
\end{equation}
\end{definition}

It can be seen that $B_{0,j}=R_{0,j}=p_{2j}$, $B_{1,j}=R_{1,j}=p_{1+2j}$, and $B_{i,0}=\Delta_i$ is the $i$-th leading principal determinant. 
For $i\geq 1$, $B_{i,j}$ is the determinant of the submatrix of $H$ formed from the first $i$ rows, and the first $i-1$ columns together with the $(i+j)$-th column. 

For instance, for $n=3$ one has
\[
\begin{array}{rcl}
H&=&\left[\begin{smallmatrix} p_1 & p_3 & 0\\ 1 & p_2 & 0\\ 0 & p_1 & p_3\end{smallmatrix}\right],\\
B_{1,0}&=&p_1,\quad B_{2,0}=p_1p_2-p_3,\quad B_{3,0}=p_3(p_1p_2-p_3).
\end{array}
\]

\begin{proposition}\label{prop:rh-bordered}
Let $A\in\mathbb R^{n\times n}$ with characteristic polynomial as in~\eqref{eq:char-poly}. Then, the Routh-Hurwitz table of $A$ is well-defined and satisfies $R_{i,j}\neq 0$ for each structural position ($i+2j\leq n$) iff $B_{i,j}\neq 0$ for each structural position. Moreover, in such cases, the entries of the Routh-Hurwitz table satisfy $R_{0,j}=B_{0,j}$ and $R_{i,j}=B_{i,j}/B_{i-1,0}$ for $i\geq 1$.~\hfill$\diamond$
\end{proposition}

\begin{myproof}
The result follows from $B_{i,j}=R_{i,j} B_{i-1,0}$, which is a consequence of Sylvester's identity~\citep{Bareiss1968} (see~(3-3) in~\cite{jeltsch1996stability}).
\end{myproof}

Let $B_{i,j}(A)$ denote the value of the bordered minor $B_{i,j}$ as a function of (the entries of) $A$. 
Notice that each $B_{i,j}(A)$ is a determinant of a submatrix of the Hurwitz matrix $H$, hence a polynomial in the entries of $H$. These entries correspond to coefficients of the characteristic polynomial $\chi_A$, and hence are polynomials in the entries of $A$. 
Consequently, $B_{i,j}(A)$ is also a polynomial in the entries of $A$. 

\begin{proposition}\label{prop:rh-bij-nonzero}
$G$ is Routh-Hurwitz iff for every structural position $(i,j)$ (with $i+2j\leq n$), there is a matrix $A\in\Sigma_G$ such that $B_{i,j}(A)\neq 0$.~\hfill$\diamond$
\end{proposition}

\begin{myproof}
Suppose that $G$ is Routh-Hurwitz. By Definition~\ref{def:scc-rh}, there is a matrix $A\in\Sigma_G$ whose Routh-Hurwitz table is well-defined and satisfies $R_{i,j}\neq 0$ for every structural position. By Proposition~\ref{prop:rh-bordered}, the same matrix satisfies $B_{i,j}(A)\neq 0$ for every such pair $(i,j)$.

Conversely, suppose that, for each structural position $(i,j)$, there is a matrix $A\in\Sigma_G$ such that $B_{i,j}(A)\neq 0$.
Then each $B_{i,j}$ is a nonzero polynomial in the free entries of $\Sigma_G$.
Since a finite product of nonzero polynomials is nonzero, the polynomial
\[
Q_G(A)=\prod_{\substack{i,j\geq 0\\i+2j\leq n}} B_{i,j}(A)
\]
is not identically zero over $\Sigma_G$.
Therefore, there is a matrix $A^\ast\in\Sigma_G$ such that $Q_G(A^\ast)\neq 0$, implying that $B_{i,j}(A^\ast)\neq 0$ simultaneously for every $(i,j)$.
By Proposition~\ref{prop:rh-bordered}, the Routh-Hurwitz table of $A^\ast$ is well-defined and satisfies $R_{i,j}\neq 0$ for every such pair.
We conclude that $G$ is Routh-Hurwitz.
\end{myproof}

We now turn our attention to the task of determining whether each $B_{i,j}$ is not identically zero.

\begin{proposition}\label{prop:computation-bordered}
Given an $n\times n$ integer matrix $A$ with entries in $\{0,\ldots,K\}$, the values $B_{i,j}(A)$ with $i+2j\leq n$ can be computed exactly in time polynomial in $n$ and $\log(K)$.
~\hfill$\diamond$
\end{proposition}

\begin{myproof}
By Proposition~\ref{prop:computation-pk} the coefficients $p_k(A)$ are integers computable in polynomial time; after forming the Hurwitz matrix $H$, each of the $\mathcal{O}(n^2)$ bordered minors can be evaluated by Bareiss Algorithm~\citep{Bareiss1968}.
\end{myproof}

\begin{lemma}\label{lem:bordered-probability-bound}
Let $G$ be a digraph on $n$ vertices and $(i,j)$ a pair of nonnegative integers with $i+2j\leq n$.
\begin{enumerate}
\item If $B_{i,j}$ is identically zero, then $\mathbb{P}_{A\sim\Sigma_G(K)}[B_{i,j}(A)=0]=1$;
\item If $B_{i,j}$ is not identically zero, then $\mathbb{P}_{A\sim\Sigma_G(K)}[B_{i,j}(A)=0]\leq n(n+1)/(2K)$.~\hfill$\diamond$
\end{enumerate}
\end{lemma}

\begin{myproof}
The first item is immediate. For the second, we notice that $B_{i,j}(A)$ is a nonzero polynomial in the entries of $A$, and we must bound its degree. 
If $i=0$, then $B_{0,j}=p_{2j}$ is either identically zero or a polynomial of degree $2j\leq n$. 
If $i\geq 1$, then each term of the determinant \eqref{eq:bordered-minor} is a signed product of entries of $H$, namely $(1,c_1),\ldots,(i,c_i)$, where $(c_1,\ldots,c_i)$ are a permutation of $\{1,\ldots,i-1\}\cup\{i+j\}$. 
Each entry $(r,c)$ of $H$ has value $p_{2c-r}$, which is either identically zero or a polynomial of degree $2c-r$. 
Hence the product of entries $(1,c_1),\ldots,(i,c_i)$ is either identically zero or a polynomial of degree 
\begin{align*}
(2c_1-1)&+(2c_2-2)+\ldots+(2c_i-i)\\
&=2(c_1+\ldots+c_i)-(1+\ldots+i)\\
&=2(1+\ldots+(i-1)+i+j)-(1+\ldots+i)\\
&=(1+\ldots+i)+2j\\
&\leq (1+\ldots+i-1)+n\\
&\leq 1+\ldots+n = n(n+1)/2,
\end{align*}
where the first inequality is justified by the constraint $i+2j\leq n$. 
Since the degree of $B_{i,j}$ is at most $n(n+1)/2$, the bound follows from the Schwartz-Zippel Lemma~\cite[Lemma A.36]{arora2009computational}.
\end{myproof}

\begin{theorem}\label{thm:rh-rp}
Determining whether a graph $G$ is Routh-Hurwitz can be decided in $\mathrm{RP}$ (probabilistic polynomial time with one-sided errors). Formally, there is a probabilistic polynomial-time algorithm that receives a digraph $G$ and: returns \true\ with probability at least $1/2$, if $G$ is Routh-Hurwitz; returns \false\ with probability $1$, if $G$ is not Routh-Hurwitz.~\hfill$\diamond$
\end{theorem}

\begin{myproof}
Consider Algorithm~\ref{alg:rh-property}, which is the exact analogue of Algorithm~\ref{alg:packability}, testing the bordered minors $B_{i,j}$ in place of the coefficients $p_k$.

\begin{algorithm}[!ht]
\caption{Randomized Routh-Hurwitz Test}
{\small
\begin{algorithmic}[1]
\Require digraph $G$ with $n$ vertices, integer $K\geq n(n+1)$, number of trials $T$
\For{$(i,j)$ such that $i+2j\leq n$}
$r_{i,j} \gets \false$
\Comment{{\color{gray}\scriptsize initialization}}
\EndFor
\For{$t=1$ to $T$}
\Comment{{\color{gray}\scriptsize repeat $T$ times}}
\State $A\gets$ random matrix in $\Sigma_{G}(K)$ \label{line:rh-start-trial}
\For{$(i,j)$ such that $i+2j\leq n$}
\If{$B_{i,j}(A)\neq 0$} $r_{i,j} \gets \true$\label{line:rh-check-zero}
\EndIf
\EndFor \label{line:rh-end-trial}
\EndFor
\State\Return{whether every $r_{i,j} = \true$} \label{line:rk-check}
\end{algorithmic}
}
\label{alg:rh-property}
\end{algorithm}

If $G$ is not Routh-Hurwitz, then by Proposition~\ref{prop:rh-bij-nonzero} some $B_{i,j}$ is identically zero, and the corresponding variable $r_{i,j}$ remains \false\ throughout all trials, so that the algorithm correctly returns \false.

Next, suppose that $G$ is Routh-Hurwitz, so that, by Proposition~\ref{prop:rh-bij-nonzero}, each $B_{i,j}$ is not identically zero.
For each trial, and each $(i,j)$, the probability that $B_{i,j}(A)=0$ (Line~\ref{line:rh-check-zero}) is bounded by $n(n+1)/(2K)\leq 1/2$, thanks to Lemma~\ref{lem:bordered-probability-bound}.
As in Theorem~\ref{thm:cycle-packing-rp}, the probability that $r_{i,j}$ remains \false\ over all $T$ trials is bounded by $2^{-T}$, and by a union bound over the $\lfloor (n+2)^2/4\rfloor$ entries, the probability that at least one $r_{i,j}$ remains \false\ is less than $\frac{(n+2)^2}{4}\cdot 2^{-T}$. Hence,  only $T \geq 2\lceil\log_2 (n+2)\rceil$ trials are required to bound this probability by $1/2$.

Taking $K=n(n+1)$, Proposition~\ref{prop:computation-bordered} guarantees that each trial (Lines~\ref{line:rh-start-trial}--\ref{line:rh-end-trial}) runs in time polynomial in $n$, so that the time complexity of Algorithm~\ref{alg:rh-property} is also polynomial.
\end{myproof}

\begin{corollary}\label{cor:scc-rh-rp}
Determining whether a graph $G$ is SCC-Routh-Hurwitz can be decided in $\mathrm{RP}$.~\hfill$\diamond$
\end{corollary}

\begin{myproof}
Similar to the proof of Corollary~\ref{cor:scc-packable-rp}.
\end{myproof}

\section{The sufficient condition is NP-complete}\label{sec:sufficient-hard}

In this section we resolve the complexity of deciding the sufficient condition (peelability), showing that it is NP-complete. 
Membership in NP is immediate, since we may take a witness for peelability as a nested sequence of cycle packings. 
We can even go one step further and merely consider permutations of the vertices of $G$ as potential witnesses, which have a more concise description. 

\begin{definition}
Let $G$ be a digraph with $n$ vertices. Let $(v_1,\ldots,v_\ell)$, with $\ell\leq n$, be a sequence of distinct vertices of $G$. 
We say that $(v_1,\ldots,v_\ell)$ is a peeling chain of $G$ if, for every $k\in\{1,\ldots,\ell\}$, the subgraph of $G$ induced by the vertices $v_1,\ldots,v_k$ has a cycle packing of size $k$.
\end{definition}

It is immediate that $G$ is peelable iff it has a complete peeling chain (of length $n$). Moreover, as already pointed out by Belabbas, we can efficiently determine whether a digraph with $k$ vertices has a cycle packing of size $k$ by reducing to the perfect matching problem.

\begin{lemma}\label{lem:full-cycle-packing}
There is a polynomial algorithm that decides, given a digraph $G$ and a subset $V'$ of the vertices of $G$, whether $G$ has a cycle packing of size $|V'|$ covering exactly the vertices in $V'$.~\hfill$\diamond$
\end{lemma}

\begin{myproof}
Let $G'$ be the subgraph of $G$ induced by the vertices of $V'$. The problem amounts to determining whether $G'$ has a full-sized cycle packing. By considering the bipartite graph $G''$ having two copies $v_{in},v_{out}$ for each $v\in V'$, and an undirected edge $\{u_{in},v_{out}\}$ for each directed edge $(u,v)$ in $G'$, it turns out that $G'$ has a full-sized cycle packing iff $G''$ has a perfect matching. Finally, we can determine whether $G''$ has a perfect matching using, for example, the Hopcroft-Karp algorithm~\citep{hopcroft1973n}.
\end{myproof}

The above lemma implies that we can efficiently decide whether a sequence of distinct vertices constitutes a peeling chain, thus proving that peelability is in NP.

\begin{lemma}\label{lem:polytime-peeling-chain}
There is a polynomial algorithm that decides, given a digraph $G$ and a sequence $(v_1,\ldots,v_\ell)$ of distinct vertices of $G$, whether $(v_1,\ldots,v_\ell)$ is a peeling chain.~\hfill$\diamond$
\end{lemma}

\begin{myproof}
Simply check if $G$ has a cycle packing of size $|V_k|$ covering exactly the vertices in $V_k$, where $V_k=\{v_1,\ldots,v_k\}$, for each $k\in\{1,\ldots,\ell\}$.
\end{myproof}

Using the above lemmas, we may design a backtracking algorithm for checking peelability (Algorithm~\ref{alg:peelability}). 


\begin{algorithm}[!ht]
\caption{Peelability Test via Memoized Backtracking}
{\small
\begin{algorithmic}[1]
\Require digraph $G=(V,E)$; subset $V'\subseteq V$; cache $\mathcal{M}$, shared across
         calls and initially empty
         \Comment{{\color{gray}\scriptsize top-level call: $V'=\{v:(v,v)\in E\}$}}
\If{$V'=V$} \Return \true \Comment{{\color{gray}\scriptsize base case}}
\EndIf
\If{$V'\in\mathcal{M}$} \Return $\mathcal{M}[V']$ \Comment{{\color{gray}\scriptsize cached verdict}}
\EndIf
\For{$v\in V\setminus V'$} \Comment{{\color{gray}\scriptsize check augmentation}}
\If {$G$ has a full cycle packing on $V'\cup\{v\}$} \label{line:full-cycle-packing}
\If {Peelability$(G,V'\cup\{v\})$} \Comment{{\color{gray}\scriptsize recursive call}}
\State \Return \true
\EndIf
\EndIf
\EndFor
\State $\mathcal{M}[V']\gets\false$;\ \Return \false
\end{algorithmic}
}
\label{alg:peelability}
\end{algorithm}


Algorithm~\ref{alg:peelability} is initialized with $V'=\{v:(v,v)\in E\}$, the set of
vertices with self-loops. At each iteration, it attempts to augment $V'$ with a new
vertex while preserving the invariant that $G$ has a full cycle packing on $V'$
(line~\ref{line:full-cycle-packing}), which reduces to perfect bipartite matching
(Lemma~\ref{lem:full-cycle-packing}). 
The self-loop initialization is without loss of
generality, since if there is a peeling chain $(v_1,\ldots,v_\ell)$ and $v_\ell$ has a self-loop, then $(v_\ell,v_1,\ldots,v_{\ell-1})$ is also a peeling chain.
When $G$ has many self-loops, this initialization removes many levels of the
backtracking search. 

If no augmentation is possible, then the algorithm
caches $V'$ in $\mathcal{M}$ before backtracking, that is, we use memoization to speed up the search. 
Whether $V'$ admits a completion to $V$
depends only on the set $V'$, and not on the order in which its vertices were added; thus, every vertex set needs to be expanded at most once and the number of recursive calls reduces from $n!$ to $n2^n$. 
If the algorithm successfully reaches $V'=V$, then the order by which the vertices were added constitutes a complete peeling chain of $G$ (of length $n$).

For the remainder of the section, we prove that deciding whether a graph is peelable is $\mathrm{NP}$-hard, by reduction from the $3$-SAT problem. 
As mentioned above, peelability is equivalent to the existence of a peeling chain $(v_1,\ldots,v_n)$ of length $n$. 


Let $\phi$ be a $3$-SAT formula with variables $x_1,\dots,x_n$ and clauses
$c_1,\dots,c_m$. We shall construct a digraph $G_\phi$ with $\mathcal O(n+m)$ vertices and $\mathcal O(n+m)$ edges, such that $\phi$ is satisfiable iff $G_\phi$ is peelable. 
The construction is split into four gadgets, which we describe sequentially. 

\begin{figure}[!ht]
    \centering
\begin{tikzpicture}[scale=.6, transform shape,node distance=1.5cm]
\begin{scope}[every node/.style={circle,thick,draw},square/.style={regular polygon,regular polygon sides=4}]
\node[label={[label distance=1mm] below:{\Large $v_0$}}] (1) at (1.25,10.) {};
\node[label={[label distance=-1mm] right:{\Large $x_1$}}] (2) at (0.,7.5) {};
\node[label={[label distance=-1mm] left:{\Large $\bar x_1$}}] (3) at (2.5,7.5) {};
\node[label={[label distance=-1mm] left:{\Large $x_2$}}] (4) at (0.,5.) {};
\node[label={[label distance=-1mm] right:{\Large $\bar x_2$}}] (5) at (2.5,5.) {};
\node[label={[label distance=-1mm] left:{\Large $x_{n-1}$}}] (6) at (0.,2.5) {};
\node[label={[label distance=-1mm] right:{\Large $\bar x_{n-1}$}}] (7) at (2.5,2.5) {};
\node[label={[label distance=0mm] left:{\Large $x_n$}}] (8) at (0.,0.) {};
\node[label={[label distance=0mm] right:{\Large $\bar x_n$}}] (9) at (2.5,0.) {};
\node[label={[label distance=1mm] above:{\Large $v_f$}}] (10) at (1.25,-2.55) {};
\end{scope}
\begin{scope}[>={Stealth[black]},
              every edge/.style={draw=black, thick}]
\path [->] (1) edge[loop above] node {} (1);
\path [->] (1) edge node {} (2);
\path [->] (1) edge node {} (3);
\path [->] (2) edge node {} (4);
\path [->] (2) edge node {} (5);
\path [->] (3) edge node {} (4);
\path [->] (3) edge node {} (5);
\path [-] (4) edge[dashed] node {} (6);
\path [-] (4) edge[dashed] node {} (7);
\path [-] (5) edge[dashed] node {} (6);
\path [-] (5) edge[dashed] node {} (7);
\path [->] (6) edge node {} (8);
\path [->] (6) edge node {} (9);
\path [->] (7) edge node {} (8);
\path [->] (7) edge node {} (9);
\path [->] (8) edge node {} (10);
\path [->] (9) edge node {} (10);
\path [->] (10) edge[bend left=85] node {} (1);
\path [->] (8) edge[bend left=70] node {} (1);
\path [->] (9) edge[bend right=70] node {} (1);
\path [->] (6) edge[bend left=60] node {} (1);
\path [->] (7) edge[bend right=60] node {} (1);
\path [->] (4) edge[bend left=50] node {} (1);
\path [->] (5) edge[bend right=50] node {} (1);
\path [->] (2) edge[bend left=30] node {} (1);
\path [->] (3) edge[bend right=30] node {} (1);
\end{scope}
\end{tikzpicture}
    \caption{Spine gadget.}
    \label{fig:chain_gadget}
\end{figure}
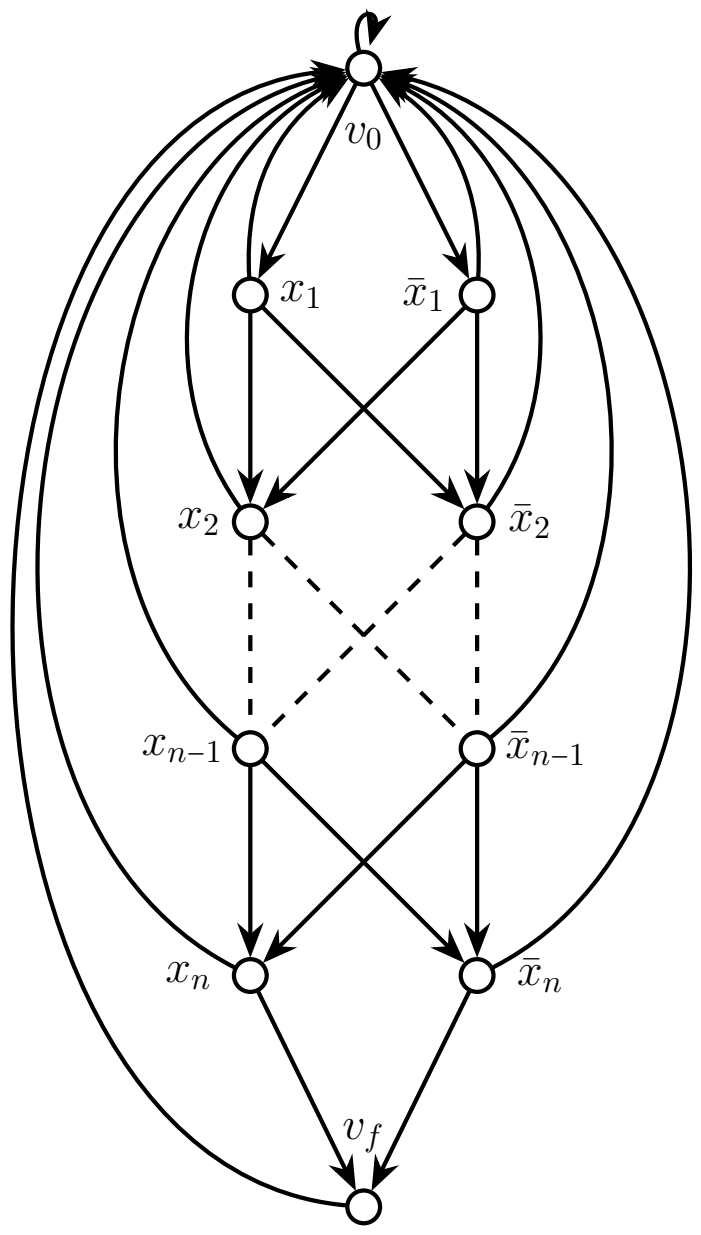

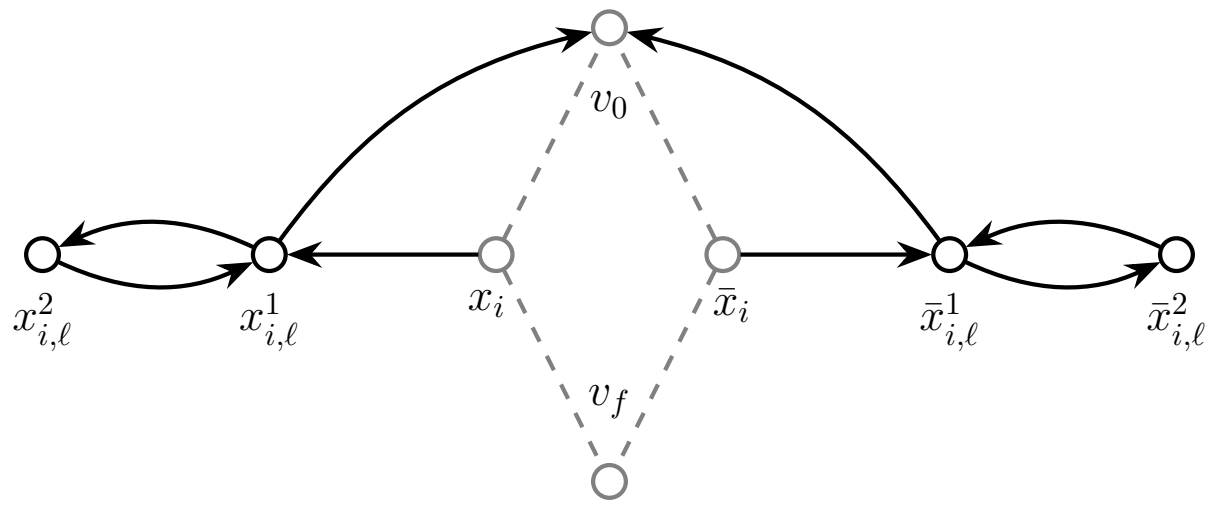
\begin{figure}[!ht]
    \centering
\begin{tikzpicture}[scale=.6, transform shape,node distance=1.5cm]
\begin{scope}[every node/.style={circle,thick,draw},square/.style={regular polygon,regular polygon sides=4}]
\node[draw=gray,label={[label distance=2mm] below:{\Large $v_0$}}] (1) at (1.25,10.) {};
\node[draw=gray,label={[label distance=-1mm,xshift=-1mm] below:{\Large $x_i$}}] (2) at (0.,7.5) {};
\node[draw=gray,label={[label distance=-1mm,xshift=1mm] below:{\Large $\bar x_i$}}] (3) at (2.5,7.5) {};
\node[label={[label distance=-1mm] below:{\Large $x_{i,\ell}^1$}}] (4) at (-2.5,7.5) {};
\node[label={[label distance=-1mm] below:{\Large $x_{i,\ell}^2$}}] (5) at (-5,7.5) {};
\node[label={[label distance=-1mm] below:{\Large $\bar x_{i,\ell}^1$}}] (6) at (5,7.5) {};
\node[label={[label distance=-1mm] below:{\Large $\bar x_{i,\ell}^2$}}] (7) at (7.5,7.5) {};
\node[draw=gray,label={[label distance=2mm] above:{\Large $v_f$}}] (8) at (1.25,5) {};
\end{scope}
\begin{scope}[>={Stealth[black]},
              every edge/.style={draw=black, thick}]
\path [->] (4) edge[bend right=25] node {} (5);
\path [->] (5) edge[bend right=25] node {} (4);
\path [->] (6) edge[bend right=25] node {} (7);
\path [->] (7) edge[bend right=25] node {} (6);
\path [->] (2) edge node {} (4);
\path [->] (3) edge node {} (6);
\path [->] (4) edge[bend left=20] node {} (1);
\path [->] (6) edge[bend right=20] node {} (1);
\path [-] (2) edge[dashed,gray] node {} (1);
\path [-] (3) edge[dashed,gray] node {} (1);
\path [-] (2) edge[dashed,gray] node {} (8);
\path [-] (3) edge[dashed,gray] node {} (8);
\end{scope}
\end{tikzpicture}
    \caption{Literal gadget.}
    \label{fig:variable_gadget}
\end{figure}

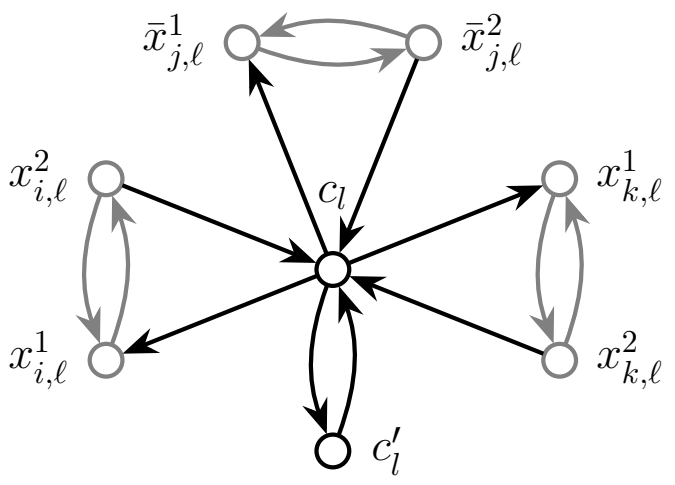
\begin{figure}[!ht]
    \centering
\begin{tikzpicture}[scale=.6, transform shape,node distance=1.5cm]
\begin{scope}[every node/.style={circle,thick,draw},square/.style={regular polygon,regular polygon sides=4}]
\node[label={[label distance=2mm] above:{\Large $c_l$}}] (1) at (0,0) {};
\node[label={[label distance=-1mm] right:{\Large $c_l'$}}] (8) at (0,-2) {};
\node[draw=gray,label={[label distance=-1mm] left:{\Large $x_{i,\ell}^1$}}] (2) at (-2.5,-1) {};
\node[draw=gray,label={[label distance=-1mm] left:{\Large $x_{i,\ell}^2$}}] (3) at (-2.5,1) {};
\node[draw=gray,label={[label distance=-1mm] left:{\Large $\bar x_{j,\ell}^1$}}] (4) at (-1,2.5) {};
\node[draw=gray,label={[label distance=-1mm] right:{\Large $\bar x_{j,\ell}^2$}}] (5) at (1,2.5) {};
\node[draw=gray,label={[label distance=-1mm] right:{\Large $x_{k,\ell}^1$}}] (6) at (2.5,1) {};
\node[draw=gray,label={[label distance=-1mm] right:{\Large $x_{k,\ell}^2$}}] (7) at (2.5,-1) {};
\end{scope}
\begin{scope}[>={Stealth[black]},
              every edge/.style={draw=black, thick}]
\path [->] (1) edge node {} (2);
\path [->] (1) edge[bend right=20] node {} (8);
\path [->] (8) edge[bend right=20] node {} (1);
\path [->] (3) edge node {} (1);
\path [->] (1) edge node {} (4);
\path [->] (5) edge node {} (1);
\path [->] (1) edge node {} (6);
\path [->] (7) edge node {} (1);
\end{scope}
\begin{scope}[>={Stealth[gray]},
              every edge/.style={draw=gray, thick}]
\path [->] (2) edge[bend right=20] node {} (3);
\path [->] (3) edge[bend right=20] node {} (2);
\path [->] (4) edge[bend right=20] node {} (5);
\path [->] (5) edge[bend right=20] node {} (4);
\path [->] (6) edge[bend right=20] node {} (7);
\path [->] (7) edge[bend right=20] node {} (6);
\end{scope}
\end{tikzpicture}
    \caption{Clause gadget. Example: $x_i\lor \bar x_j\lor x_k$.}
    \label{fig:clause_gadget}
\end{figure}

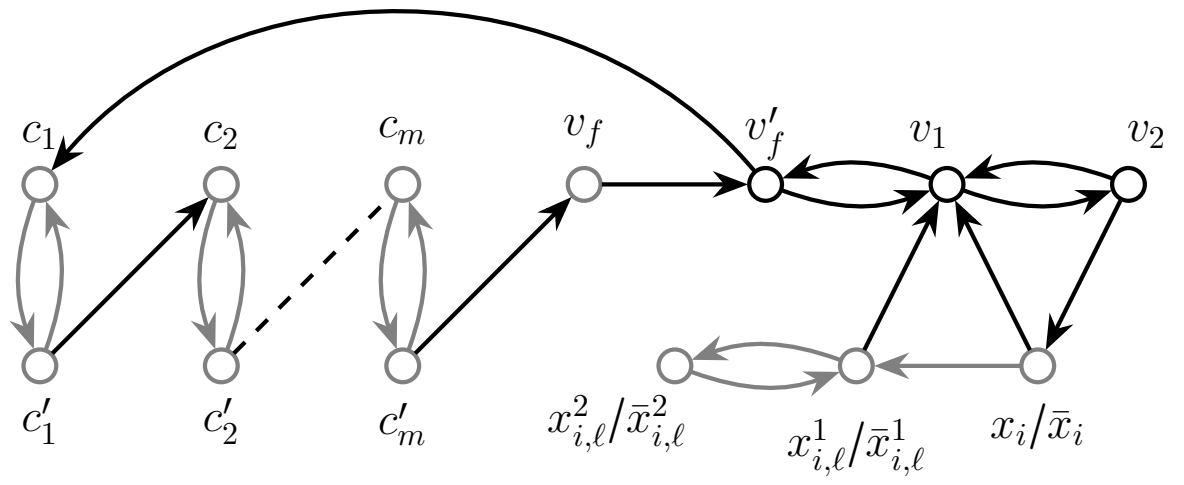
\begin{figure}[!ht]
    \centering
\begin{tikzpicture}[scale=.6, transform shape,node distance=1.5cm]
\begin{scope}[every node/.style={circle,thick,draw},square/.style={regular polygon,regular polygon sides=4}]
\node[draw=gray,label={[label distance=-1mm] above:{\Large $c_1$}}] (1) at (0,0) {};
\node[draw=gray,label={[label distance=-1mm] below:{\Large $c_1'$}}] (1b) at (0,-2) {};
\node[draw=gray,label={[label distance=-1mm] above:{\Large $c_2$}}] (2) at (2,0) {};
\node[draw=gray,label={[label distance=-1mm] below:{\Large $c_2'$}}] (2b) at (2,-2) {};
\node[draw=gray,label={[label distance=-1mm] above:{\Large $c_m$}}] (3) at (4,0) {};
\node[draw=gray,label={[label distance=-1mm] below:{\Large $c_m'$}}] (3b) at (4,-2) {};
\node[draw=gray,label={[label distance=-1mm] above:{\Large $v_f$}}](5) at (6,0) {};
\node[label={[label distance=-2mm] above:{\Large $v_f'$}}]  (6) at (8,0) {};
\node[label={[xshift=-2mm,yshift=-1mm] above:{\Large $v_1$}}] (7) at (10,0) {};
\node[label={[xshift=2mm,yshift=-1mm] above:{\Large $v_2$}}] (8) at (12,0) {};
\node[draw=gray,label={[label distance=-3mm] below:{\Large $x_{i}/\bar x_{i}$}}] (9) at (11,-2) {};
\node[draw=gray,label={[label distance=-3mm] below:{\Large $x_{i,\ell}^1/\bar x_{i,\ell}^1$}}] (10) at (9,-2) {};
\node[draw=gray,label={[label distance=-3mm] below left:{\Large $x_{i,\ell}^2/\bar x_{i,\ell}^2$}}] (11) at (7,-2) {};
\end{scope}
\begin{scope}[>={Stealth[black]},
              every edge/.style={draw=black, thick}]
\path [->] (1b) edge node {} (2);
\path [-] (2b) edge[dashed] node {} (3);
\path [->] (3b) edge node {} (5);
\path [->] (5) edge node {} (6);
\path [->] (6) edge[bend right=50] node {} (1);
\path [->] (6) edge[bend right=20] node {} (7);
\path [->] (7) edge[bend right=20] node {} (6);
\path [->] (7) edge[bend right=20] node {} (8);
\path [->] (8) edge[bend right=20] node {} (7);
\path [->] (9) edge node {} (7);
\path [->] (8) edge node {} (9);
\path [->] (10) edge node {} (7);
\end{scope}
\begin{scope}[>={Stealth[gray]},
              every edge/.style={draw=gray, thick}]
\path [->] (9) edge node {} (10);
\path [->] (10) edge[bend right=20] node {} (11);
\path [->] (11) edge[bend right=20] node {} (10);
\path [->] (1) edge[bend right=20] node {} (1b);
\path [->] (1b) edge[bend right=20] node {} (1);
\path [->] (2) edge[bend right=20] node {} (2b);
\path [->] (2b) edge[bend right=20] node {} (2);
\path [->] (3) edge[bend right=20] node {} (3b);
\path [->] (3b) edge[bend right=20] node {} (3);
\end{scope}
\end{tikzpicture}
    \caption{Validation gadget.}
    \label{fig:final_gadget}
\end{figure}

The intuition is as follows. 
A peeling chain is forced to descend a \emph{spine}, and each descent selects, for every variable, one of two \emph{choice vertices}: this encodes a truth assignment, one choice per variable. 
Each choice vertex is attached to \emph{literal
vertex pairs}, one for each occurrence of the literal in a clause, which can be added to the chain after the corresponding choice vertex.
Each \emph{clause vertex} can be covered through a small cycle via one of its literal pairs, so a clause vertex can only be added to the chain when one of its literals
has been chosen true. 
A validation gadget forces the peeling chain to cover all clause
vertices before it can cover a \emph{validation vertex} and a \emph{cleanup vertex pair}, which can be used to cover all remaining uncovered vertices. 
Thus a peeling chain exists if and only if some
assignment makes every clause true, that is, if and only if $\phi$ is satisfiable. 
The four gadgets, shown in
Figures~\ref{fig:chain_gadget}--\ref{fig:final_gadget}, realize these four roles: spine, literal, clause, and validation.

\emph{Spine gadget} (Figure~\ref{fig:chain_gadget}). 
We include in $G_\phi$ a root vertex $v_0$, which is the only vertex with a self-loop, as well as a final vertex $v_f$. 
For each variable $x_i$, we include two choice vertices $x_i$ and $\bar x_i$. 
The spine gadget includes the initial edges
$v_0\to x_1$, $v_0\to\bar x_1$; 
the crossing edges
$x_i,\bar x_i\to x_{i+1},\bar x_{i+1}$ between consecutive levels; 
the final edges $x_n \to v_f$, $\bar x_n\to v_f$;
and the returning edges $x_i\to v_0$, $\bar x_i\to v_0$, $v_f\to v_0$. 

Notice that every peeling chain must begin at $v_0$, since it is the only vertex with a self-loop.
Moreover, every cycle in the spine gadget must include $v_0$, so that every cycle packing in the spine gadget must consist of a single cycle. 
Since the edges leaving level $i$ go to level $i+1$ or back to $v_0$, we conclude that any peeling chain adds exactly one choice vertex per level.
Therefore, there are exactly $2^n$ peeling chains of length $n+2$ in the spine gadget; these are in bijection with the truth assignments of $\phi$, where intuitively, entering vertex $x_i$ corresponds to $x_i=\true$ and entering vertex $\bar x_i$ corresponds to $x_i=\false$.

\emph{Literal gadgets} (Figure~\ref{fig:variable_gadget}). For each choice vertex, say $x_i$, 
and for each occurrence of that literal in a clause, say $c_\ell$, 
we include a pair of literal vertices $x_{i,\ell}^1,x_{i,\ell}^2$ joined by a $2$-cycle. 
Each literal gadget also includes the entering edge $x_{i}\to x_{i,\ell}^1$ and the returning edge $x_{i,\ell}^1\to v_0$. 
Similarly, for a negative literal $\bar x_i$ occurring in clause $c_\ell$, we have $\bar x_{i}\to \bar x_{i,\ell}^1 \to v_0$.


Consider one of the $2^n$ peeling chains of length $n+2$ in the spine gadget. 
We may augment the chain by including each literal vertex pair associated to each occurrence of the choice vertex in each clause. 
Formally, after adding $x_i$ but before moving to level $i+1$, we may add $x_{i,\ell}^1$ (via the cycle $v_0\to\ldots\to x_i\to x_{i,\ell}^1 \to v_0$) followed by $x_{i,\ell}^2$ (via the 2-cycle on $\{x_{i,\ell}^1,x_{i,\ell}^2\}$, plus the original cycle in the spine gadget) and repeat for every occurrence of $x_i$. 
For negative literals the process is similar.


\emph{Clause gadgets} (Figure~\ref{fig:clause_gadget}). 
For each clause $c_\ell$, we include a pair of
clause vertices labelled $c_\ell,c_\ell'$, joined by a 2-cycle. 
For each literal of $c_\ell$, say $x_i$, each clause gadget also includes the edges $x_{i,\ell}^2\to c_\ell$ and $c_\ell\to x_{i,\ell}^1$.

Consider a peeling chain on the spine and literal gadgets, as before. We may augment the chain by including the clause vertex pairs associated with the satisfied clauses. 
Formally, if the peeling chain contains the literal vertex pair $\{x_{i,\ell}^1,x_{i,\ell}^2\}$, and $x_i\in c_\ell$, we may add $c_\ell$ to the chain by replacing the 2-cycle with a 3-cycle on $\{x_{i,\ell}^1,x_{i,\ell}^2,c_\ell\}$, followed by $c_\ell'$ via its 2-cycle. For negative literals the process is similar.



\emph{Validation gadget} (Figure~\ref{fig:final_gadget}). We include a validation vertex $v_f'$, which is chained with the final vertex $v_f$ and the clause vertices into a cycle $c_1\to c_1'\to\cdots\to c_m\to c_m'\to v_f\to v_f'\to c_1$. 
Furthermore, we include a pair of cleanup vertices $v_1,v_2$ joined by a 2-cycle. 
The validation gadget also includes the 2-cycle on $\{v_f', v_1\}$, as well as the cleanup edges $v_2\to x_i$, $x_i \to v_1$ for each choice vertex $x_i$, the edge $x_{i,\ell}^1 \to v_1$, for each literal vertex $x_{i,\ell}^1$, and similarly for the negative vertices $\bar x_i, \bar x_{i,\ell}^1$.

Consider a peeling chain on the spine, literal and clause gadgets that covers all clause vertices, as before. 
We may augment the chain to include $v_f'$, using the cycle that traverses $v_f$ and all the clauses; then include $v_1$ using the 2-cycle $\{v_f',v_1\}$; then include $v_2$ using its 2-cycle. 
From this point on, we may use the cleanup edges to absorb all unused choice vertices and literal vertices. 
Formally, let $u_1,\ldots u_n$ be the unused choice vertices; each $u_i$ is either $x_i$ or $\bar x_i$. 
At level $i$, we may add the unused choice vertex $u_i$ via the cycle $v_1\to v_2\to u_1\to \cdots \to u_i\to v_1$; then, for each unused literal vertex pair $\{u_{i,\ell}^1,u_{i,\ell}^2\}$, we first add $u_{i,\ell}^1$ via the similar cycle that ends in $\ldots\to u_i \to u_{i,\ell}^1 \to v_1$; finally, we add $u_{i,\ell}^2$ via its 2-cycle. 

\begin{lemma}\label{lem:sat-peelable}
If $\phi$ is satisfiable, then $G_\phi$ is peelable.~\hfill$\diamond$
\end{lemma}

\begin{myproof}
Immediate from the above discussion.
\end{myproof}

We now consider the converse of Lemma~\ref{lem:sat-peelable}, and show that if $G_\phi$ is peelable, then $\phi$ is satisfiable. 
We need to argue that any peeling chain must, essentially, follow the order outlined above.

\begin{lemma}\label{lem:vf-before-vi}
Any peeling chain on $G_\phi$ can only include the cleanup vertices $v_1,v_2$ after including the validation vertex $v_f'$.~\hfill$\diamond$
\end{lemma}

\begin{myproof}
Let $C$ be a peeling chain on $G_\phi$. Since $v_1 \to v_2$ is the only edge entering $v_2$, every cycle containing $v_2$ must also include $v_1$. 
Therefore, $v_2$ may only be added after $v_1$. 
Now suppose $v_1$ is added to $C$. 
By the previous observation, the cycle used to add $v_1$ cannot use edge $v_1 \to v_2$, and thus it must use the only other edge exiting $v_1$, which is $v_1\to v_f'$. We conclude that $v_f'$ must have been added before $v_1$.
\end{myproof}

\begin{lemma}\label{lem:none-both-literals}
Let $C$ be a peeling chain on $G_\phi$ that includes $v_f$ but does not include $v_f'$. Then, for each positive literal vertex pair $x_{i,\ell}^1,x_{i,\ell}^2$, $C$ must either include none or both of the vertices; and similarly for each negative literal vertex pair $\bar x_{i,\ell}^1,\bar x_{i,\ell}^2$.~\hfill$\diamond$
\end{lemma}

\begin{myproof}
Consider a positive literal vertex pair $x_{i,\ell}^1,x_{i,\ell}^2$. The only edge entering $x_{i,\ell}^2$ is $x_{i,\ell}^1\to x_{i,\ell}^2$; therefore, if $x_{i,\ell}^2$ is in $C$, then $x_{i,\ell}^1$ must have been previously added to $C$.
Now suppose that $x_{i,\ell}^1$ is in $C$.
Consider a cycle packing covering all vertices in $C$; this cycle packing must include some edge exiting $x_{i,\ell}^1$. Since $C$ includes $v_f$ but does not include $v_f'$, the cycle packing must use edge $v_f\to v_0$, and cannot therefore use edge $x_{i,\ell}^1\to v_0$. Moreover, by Lemma~\ref{lem:vf-before-vi}, $C$ also does not include $v_1$, and cannot therefore use edge $x_{i,\ell}^1\to v_1$. The only remaining possibility is for the cycle packing to use edge $x_{i,\ell}^1\to x_{i,\ell}^2$, implying that $x_{i,\ell}^2$ is also in $C$. For negative literal vertex pairs, the argument is similar.
\end{myproof}

\begin{lemma}\label{lem:cl-before-vf}
Any peeling chain on $G_\phi$ can only include the validation vertex $v_f'$ after including the final vertex $v_f$ as well as all clause vertices $c_1,c_1',\ldots,c_m,c_m'$.~\hfill$\diamond$
\end{lemma}

\begin{myproof}
Let $C$ be a peeling chain on $G_\phi$ that includes $v_f'$ and let $C'$ be the prefix of $C$ immediately before $v_f'$ is added. 
By Lemma~\ref{lem:vf-before-vi}, $C'$ cannot include vertices $v_1, v_2$. Thus, the cycle used to add $v_f'$ must have used the only other edges adjacent to $v_f'$, which are $v_f'\to c_1$ and $v_f\to v_f'$. 
In particular, $v_f$ must belong to $C'$.
Suppose that the cycle used to add $v_f'$ is not $v_f'\to c_1\to c_1' \to \cdots \to c_m \to c_m' \to v_f\to v_f'$, and let $c_\ell$ be the first clause vertex at which point the cycle deviates (it cannot first deviate at $c_\ell'$ since the only other edge leaving $c_\ell'$ is $c_\ell'\to c_\ell$, which has already been visited). 
The cycle must then use edge $c_\ell\to u_\ell^1$, where $u_\ell^1$ is a literal vertex associated with $c_\ell$ (either $x_{i,\ell}^1$ or $\bar x_{i,\ell}^1$ for some $i$). 
By Lemma~\ref{lem:none-both-literals}, the literal vertex $u_\ell^2$ must also appear in $C'$, but its only two outgoing edges lead back to $c_\ell$ and $u_\ell^1$, which are already in a cycle; a contradiction.   
We conclude that all clause vertices must belong to $C'$.
\end{myproof}

\begin{lemma}\label{lem:peelable-sat}
If $G_\phi$ is peelable, then $\phi$ is satisfiable.~\hfill$\diamond$
\end{lemma}

\begin{myproof}
Suppose that $G_\phi$ is peelable and let $C$ be a complete peeling chain on $G_\phi$. 
Consider the prefix $C'$ of $C$ immediately before the validation vertex $v_f'$ is added; by Lemmas~\ref{lem:vf-before-vi} and~\ref{lem:cl-before-vf}, $v_f$ and all $c_1,c_1',\ldots,c_m,c_m'$ have been added, but $v_1,v_2$ have not yet been added. 
Since the only edge entering $c_\ell'$ is $c_\ell\to c_\ell'$, the vertex $c_\ell$ must be added before $c_\ell'$. 
The only other exiting edges from $c_\ell$ are of the form $c_\ell\to u_\ell^1$, for $u_\ell$ an occurrence (either $x_{i,\ell}$ or $\bar x_{i,\ell}$) of a literal $u$ (either $x_i$ or $\bar x_i$) in $c_\ell$. 
Hence, at least one such literal vertex $u_\ell^1$ must have been added before $c_\ell$. 
By Lemma~\ref{lem:none-both-literals}, the literal vertex $u_\ell^2$ (either $x_{i,\ell}^2$ or $\bar x_{i,\ell}^2$) must also be present in $C'$; and moreover, $u_\ell^1$ was added before $u_\ell^2$. 
At the step in which $u_\ell^1$ was added, neither $u_\ell^2$ nor $c_\ell$ was present. 
The only remaining edge entering $u_\ell^1$ is $u\to u_\ell^1$, so the choice vertex $u$ must have been added before $u_\ell^1$. 
In summary, for each clause $c_\ell$ there exists a literal $u\in c_\ell$ such that the choice vertex $u$ is in the chain~$C'$. 

Next, we argue that $C'$ cannot contain both a choice vertex $x_i$ and its negation $\bar x_i$. Since the vertices $v_1,v_2$ are not in $C'$, any cycle that traverses $x_i$ must pass through the root vertex $v_0$; and similarly for $\bar x_i$. 
Moreover, any path between $x_i$ and $\bar x_i$, in either direction, must pass through $v_0$ as well. 
Since $v_0$ can only be used once in a cycle packing, we cannot cover both $x_i$ and $\bar x_i$. 
We conclude that $C'$ contains at most one of the choice vertices $x_i,\bar x_i$, and that it contains at least one choice vertex associated to a literal in each clause. Thus, the choice vertices in $C'$ can be used to define a satisfying assignment for $\phi$.
\end{myproof}

\begin{theorem}\label{thm:np}
Determining whether a digraph $G$ is peelable is $\mathrm{NP}$-complete.~\hfill$\diamond$
\end{theorem}

\begin{myproof}
Membership in $\mathrm{NP}$ follows from Lemma~\ref{lem:polytime-peeling-chain}. 
$\mathrm{NP}$-hardness follows from the reduction described above, since Lemmas~\ref{lem:sat-peelable} and~\ref{lem:peelable-sat} show that $\phi$ is satisfiable iff $G_\phi$ is peelable. Moreover, constructing $G_\phi$ from $\phi$ can be done in polynomial time; in particular, if $\phi$ has $n$ variables, $m$ clauses, and at most $3$ literals per clause, then $G_\phi$ has at most $5+2n+8m$ vertices and at most $8+10n+24m$ edges (linear size).
%
%
%
\end{myproof}

\section{Experimental results}\label{sec:experiments}

We implemented the conditions in Section~\ref{sec:strict} as well as a stability test,
and ran three experiments: an exhaustive enumeration for small $n$, a
phase-transition study for random graphs, and a scaling study of the running
times. 
The sufficient condition (peelability) is checked according to Algorithm~\ref{alg:peelability}. 
The necessary conditions (SCC-packability and SCC-RH) are checked probabilistically according to Algorithms~\ref{alg:packability} and~\ref{alg:rh-property}, sampling integers between $1$ and $K=2^{30}$, and taking $T=10$ trials. Since the largest graphs considered have $n=30$ vertices, the probability of a false negative is less than $\frac{32^2}{4}\times\left(\frac{30\times 31}{2^{31}}\right)^{10}< 10^{-61}$, according to the proof of Theorem~\ref{thm:rh-rp}~\footnote{Counting over all our experiments, we ran the probabilistic algorithms in fewer than $10^{8}$ digraphs, so the overall probability of a false negative is below $10^{-53}$.}. 


For the graphs in the gap between SCC-RH and peelability, and for $n\leq 5$, we further check for stability. 
We first perform a numerical search, setting random values between $-4$ and $4$ on the free entries, and checking if all eigenvalues have real part below $-10^{-9}$ (this threshold is used to remove floating-point inaccuracies that return near-zeros); a successful check provides a Hurwitz matrix for the graph. 
If the numerical search fails, we resort to the z3 SMT solver~\citep{de2008z3}, by writing the characteristic polynomial coefficients $p_k$ and Hurwitz determinants $\Delta_k$ as expressions on the entries $a_{ij}$ (which are set to unrestricted reals), and searching for instances satisfying all $p_k,\Delta_k>0$. All calls to z3 returned \texttt{sat} or \texttt{unsat}, without timeouts. 

\subsection{Exhaustive enumeration}

For $2\leq n\leq 5$ we enumerated all $2^{n^2}$ digraphs on $n$ vertices and counted
how many satisfy each property, reported in Table~\ref{tab:exhaustive}.
The nested inclusions
$\text{peelable}\subseteq\text{stable}\subseteq\text{SCC-RH}\subseteq
\text{SCC-packable}$ hold exactly in every row, as they must.

\begin{table}[!ht]
\centering
\caption{Number of digraphs on $n$ vertices satisfying each property, out of all
$2^{n^2}$ digraphs.}
\label{tab:exhaustive}
\resizebox{.75\columnwidth}{!}{%
\begin{tabular}{@{}rrrrrr@{}}
\toprule
$n$ & total $2^{n^2}$ & SCC-pack. & SCC-RH & stable & peelable \\
\midrule
2 & 16          & 6          & 6          & 6          & 6 \\
3 & 512         & 205        & 205        & 205        & 199 \\
4 & 65\,536     & 33\,428    & 33\,356    & 33\,356    & 31\,328 \\
5 & 33\,554\,432 & 21\,817\,221 & 21\,793\,941 & 21\,787\,581 & 20\,445\,031 \\
\bottomrule
\end{tabular}
}
\end{table}

At $n=2$, the four properties
coincide. 
At
$n=3$, the gap between stability and peelability emerges: there are~$6$ digraphs that are stable but not peelable, one being the witness of Proposition~\ref{prop:notnec},
and this gap grows
to about $4.00\%$ of all graphs at $n=5$. 
The gap between SCC-packability and SCC-RH opens next, at $n=4$, with 72 digraphs in the gap, including the witness of
Proposition~\ref{prop:packable-not-rh}; this gap amounts to about $0.07\%$ of all graphs at $n=5$. 
The last gap to open, at $n=5$, is the one between the
SCC-RH and stability: there are 6360 graphs which are SCC-RH but unstable, one being the witness of
Proposition~\ref{prop:rh-not-sufficient}; the gap amounts to about $0.02\%$ of all graphs at $n=5$. 
From these results, we conclude that the SCC-RH condition is a decent approximation of stability, at least for small graphs.

\subsection{Phase transition}

For $n=5$, we group the digraphs by their number of edges $m$ (including self-loops), and computed the proportion of graphs which are SCC-RH, stable, and peelable. The three conditions undergo a phase transition
at $m\approx 12$. For sufficiently sparse ($m<7$)
or dense ($m>17$) graphs, all three conditions agree.

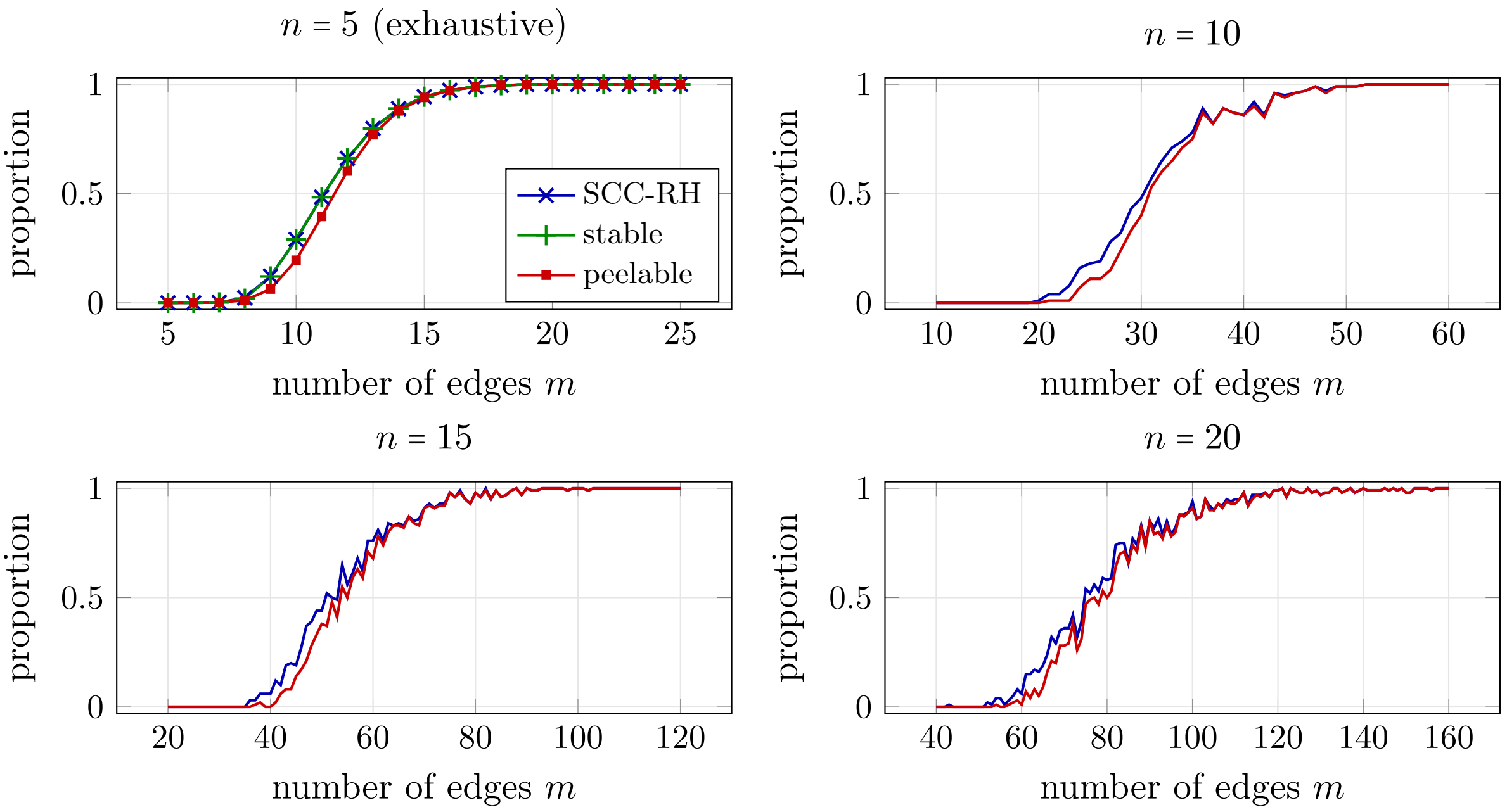
\begin{figure*}[!tbp]\centering
\begin{tikzpicture}
\begin{groupplot}[
  group style={group size=2 by 2, horizontal sep=1.6cm, vertical sep=1.8cm},
  width=8cm, height=4cm,ymin=-0.03, ymax=1.03,
  xlabel={number of edges $m$}, ylabel={proportion},
  grid=both, grid style={gray!20}, tick label style={font=\small},
  legend style={font=\footnotesize, at={(0.98,0.03)}, anchor=south east},
  legend cell align=left]
\nextgroupplot[title={$n=5$ (exhaustive)}]
\addplot[mark=x, mark size=3pt, thick, color=blue!70!black] coordinates {(5,0.0000) (6,0.0002) (7,0.0021) (8,0.0229) (9,0.1218) (10,0.2902) (11,0.4838) (12,0.6612) (13,0.7979) (14,0.8890) (15,0.9432) (16,0.9729) (17,0.9883) (18,0.9957) (19,0.9988) (20,0.9998) (21,1.0000) (22,1.0000) (23,1.0000) (24,1.0000) (25,1.0000)};
\addplot[mark=+, mark size=3pt, thick, color=green!55!black] coordinates {(5,0.0000) (6,0.0002) (7,0.0021) (8,0.0192) (9,0.1208) (10,0.2901) (11,0.4838) (12,0.6612) (13,0.7979) (14,0.8890) (15,0.9432) (16,0.9729) (17,0.9883) (18,0.9957) (19,0.9988) (20,0.9998) (21,1.0000) (22,1.0000) (23,1.0000) (24,1.0000) (25,1.0000)};
\addplot[mark=square*, mark size=1pt, thick, color=red!80!black] coordinates {(5,0.0000) (6,0.0002) (7,0.0018) (8,0.0116) (9,0.0629) (10,0.1951) (11,0.3951) (12,0.6032) (13,0.7693) (14,0.8783) (15,0.9401) (16,0.9723) (17,0.9882) (18,0.9957) (19,0.9988) (20,0.9998) (21,1.0000) (22,1.0000) (23,1.0000) (24,1.0000) (25,1.0000)};
\legend{SCC-RH, stable, peelable}
\nextgroupplot[title={$n=10$}]
\addplot[thick, color=blue!70!black] coordinates {
(10,0.0000) (11,0.0000) (12,0.0000) (13,0.0000) (14,0.0000) (15,0.0000) (16,0.0000) (17,0.0000) (18,0.0000) (19,0.0000) (20,0.0100) (21,0.0400) (22,0.0400) (23,0.0800) (24,0.1600) (25,0.1800) (26,0.1900) (27,0.2800) (28,0.3200) (29,0.4300) (30,0.4800) (31,0.5700) (32,0.6500) (33,0.7100) (34,0.7400) (35,0.7800) (36,0.8900) (37,0.8200) (38,0.8900) (39,0.8700) (40,0.8600) (41,0.9200) (42,0.8600) (43,0.9600) (44,0.9500) (45,0.9600) (46,0.9700) (47,0.9900) (48,0.9700) (49,0.9900) (50,0.9900) (51,0.9900) (52,1.0000) (53,1.0000) (54,1.0000) (55,1.0000) (56,1.0000) (57,1.0000) (58,1.0000) (59,1.0000) (60,1.0000) 
};
\addplot[thick, color=red!80!black] coordinates {
(10,0.0000) (11,0.0000) (12,0.0000) (13,0.0000) (14,0.0000) (15,0.0000) (16,0.0000) (17,0.0000) (18,0.0000) (19,0.0000) (20,0.0000) (21,0.0100) (22,0.0100) (23,0.0100) (24,0.0700) (25,0.1100) (26,0.1100) (27,0.1500) (28,0.2400) (29,0.3300) (30,0.4000) (31,0.5300) (32,0.6000) (33,0.6500) (34,0.7100) (35,0.7500) (36,0.8700) (37,0.8200) (38,0.8900) (39,0.8700) (40,0.8600) (41,0.9000) (42,0.8500) (43,0.9600) (44,0.9400) (45,0.9600) (46,0.9700) (47,0.9900) (48,0.9600) (49,0.9900) (50,0.9900) (51,0.9900) (52,1.0000) (53,1.0000) (54,1.0000) (55,1.0000) (56,1.0000) (57,1.0000) (58,1.0000) (59,1.0000) (60,1.0000) 
};
\nextgroupplot[title={$n=15$}]
\addplot[color=blue!70!black, thick] coordinates {
(20,0.0000) (21,0.0000) (22,0.0000) (23,0.0000) (24,0.0000) (25,0.0000) (26,0.0000) (27,0.0000) (28,0.0000) (29,0.0000) (30,0.0000) (31,0.0000) (32,0.0000) (33,0.0000) (34,0.0000) (35,0.0000) (36,0.0300) (37,0.0300) (38,0.0600) (39,0.0600) (40,0.0600) (41,0.1200) (42,0.1000) (43,0.1900) (44,0.2000) (45,0.1900) (46,0.2700) (47,0.3700) (48,0.3900) (49,0.4400) (50,0.4400) (51,0.5200) (52,0.5000) (53,0.4900) (54,0.6500) (55,0.5600) (56,0.6100) (57,0.6800) (58,0.6200) (59,0.7600) (60,0.7600) (61,0.8100) (62,0.7600) (63,0.8400) (64,0.8300) (65,0.8400) (66,0.8300) (67,0.8700) (68,0.8500) (69,0.8600) (70,0.9100) (71,0.9300) (72,0.9100) (73,0.9300) (74,0.9300) (75,0.9800) (76,0.9600) (77,0.9900) (78,0.9500) (79,0.9300) (80,0.9800) (81,0.9600) (82,1.0000) (83,0.9500) (84,0.9900) (85,0.9600) (86,0.9700) (87,0.9900) (88,1.0000) (89,0.9700) (90,1.0000) (91,0.9900) (92,0.9900) (93,1.0000) (94,1.0000) (95,1.0000) (96,1.0000) (97,1.0000) (98,0.9900) (99,1.0000) (100,1.0000) (101,1.0000) (102,0.9900) (103,1.0000) (104,1.0000) (105,1.0000) (106,1.0000) (107,1.0000) (108,1.0000) (109,1.0000) (110,1.0000) (111,1.0000) (112,1.0000) (113,1.0000) (114,1.0000) (115,1.0000) (116,1.0000) (117,1.0000) (118,1.0000) (119,1.0000) (120,1.0000) 
};
\addplot[thick, color=red!80!black] coordinates {
(20,0.0000) (21,0.0000) (22,0.0000) (23,0.0000) (24,0.0000) (25,0.0000) (26,0.0000) (27,0.0000) (28,0.0000) (29,0.0000) (30,0.0000) (31,0.0000) (32,0.0000) (33,0.0000) (34,0.0000) (35,0.0000) (36,0.0000) (37,0.0100) (38,0.0200) (39,0.0000) (40,0.0000) (41,0.0200) (42,0.0600) (43,0.0800) (44,0.0800) (45,0.1400) (46,0.1700) (47,0.2100) (48,0.2800) (49,0.3300) (50,0.3800) (51,0.3700) (52,0.4800) (53,0.4100) (54,0.5500) (55,0.5000) (56,0.5900) (57,0.6300) (58,0.5900) (59,0.7100) (60,0.6800) (61,0.7800) (62,0.7400) (63,0.8000) (64,0.8300) (65,0.8300) (66,0.8200) (67,0.8700) (68,0.8400) (69,0.8300) (70,0.9100) (71,0.9200) (72,0.9100) (73,0.9200) (74,0.9200) (75,0.9800) (76,0.9600) (77,0.9800) (78,0.9500) (79,0.9300) (80,0.9800) (81,0.9600) (82,0.9900) (83,0.9500) (84,0.9900) (85,0.9600) (86,0.9700) (87,0.9900) (88,1.0000) (89,0.9700) (90,1.0000) (91,0.9900) (92,0.9900) (93,1.0000) (94,1.0000) (95,1.0000) (96,1.0000) (97,1.0000) (98,0.9900) (99,1.0000) (100,1.0000) (101,1.0000) (102,0.9900) (103,1.0000) (104,1.0000) (105,1.0000) (106,1.0000) (107,1.0000) (108,1.0000) (109,1.0000) (110,1.0000) (111,1.0000) (112,1.0000) (113,1.0000) (114,1.0000) (115,1.0000) (116,1.0000) (117,1.0000) (118,1.0000) (119,1.0000) (120,1.0000) 
};
\nextgroupplot[title={$n=20$}]
\addplot[thick, color=blue!70!black] coordinates {
(40,0.0000) (41,0.0000) (42,0.0000) (43,0.0100) (44,0.0000) (45,0.0000) (46,0.0000) (47,0.0000) (48,0.0000) (49,0.0000) (50,0.0000) (51,0.0000) (52,0.0200) (53,0.0100) (54,0.0400) (55,0.0400) (56,0.0100) (57,0.0300) (58,0.0500) (59,0.0800) (60,0.0600) (61,0.1500) (62,0.1500) (63,0.1700) (64,0.1600) (65,0.1900) (66,0.2400) (67,0.3200) (68,0.2900) (69,0.3500) (70,0.3600) (71,0.3600) (72,0.4200) (73,0.3200) (74,0.3900) (75,0.5400) (76,0.5200) (77,0.5600) (78,0.5300) (79,0.5900) (80,0.5800) (81,0.5900) (82,0.7400) (83,0.7500) (84,0.7500) (85,0.6600) (86,0.7700) (87,0.7400) (88,0.8300) (89,0.7500) (90,0.8500) (91,0.8200) (92,0.8600) (93,0.7900) (94,0.8500) (95,0.7900) (96,0.8200) (97,0.8800) (98,0.8800) (99,0.8900) (100,0.9400) (101,0.8600) (102,0.8700) (103,0.9500) (104,0.9200) (105,0.9000) (106,0.9300) (107,0.9200) (108,0.9500) (109,0.9400) (110,0.9500) (111,0.9500) (112,0.9800) (113,0.9200) (114,0.9700) (115,0.9700) (116,0.9700) (117,0.9800) (118,0.9600) (119,0.9900) (120,0.9900) (121,1.0000) (122,0.9600) (123,1.0000) (124,0.9900) (125,0.9800) (126,0.9800) (127,1.0000) (128,0.9800) (129,0.9900) (130,0.9700) (131,0.9800) (132,0.9800) (133,1.0000) (134,1.0000) (135,0.9800) (136,0.9900) (137,1.0000) (138,0.9800) (139,0.9900) (140,1.0000) (141,0.9900) (142,0.9900) (143,0.9900) (144,0.9900) (145,1.0000) (146,0.9900) (147,1.0000) (148,0.9900) (149,1.0000) (150,0.9800) (151,0.9800) (152,1.0000) (153,1.0000) (154,1.0000) (155,1.0000) (156,0.9900) (157,1.0000) (158,1.0000) (159,1.0000) (160,1.0000) 
};
\addplot[thick, color=red!80!black] coordinates {
(40,0.0000) (41,0.0000) (42,0.0000) (43,0.0000) (44,0.0000) (45,0.0000) (46,0.0000) (47,0.0000) (48,0.0000) (49,0.0000) (50,0.0000) (51,0.0000) (52,0.0000) (53,0.0000) (54,0.0100) (55,0.0000) (56,0.0000) (57,0.0100) (58,0.0200) (59,0.0300) (60,0.0100) (61,0.0700) (62,0.0400) (63,0.0800) (64,0.0500) (65,0.0900) (66,0.1600) (67,0.2100) (68,0.2000) (69,0.2800) (70,0.2800) (71,0.2900) (72,0.3800) (73,0.2600) (74,0.3100) (75,0.4700) (76,0.4900) (77,0.5000) (78,0.4700) (79,0.5300) (80,0.5000) (81,0.5300) (82,0.6400) (83,0.7000) (84,0.7100) (85,0.6600) (86,0.7400) (87,0.7100) (88,0.8200) (89,0.7300) (90,0.8500) (91,0.7900) (92,0.8000) (93,0.7700) (94,0.8300) (95,0.7800) (96,0.8000) (97,0.8800) (98,0.8700) (99,0.8900) (100,0.9100) (101,0.8600) (102,0.8700) (103,0.9500) (104,0.9000) (105,0.9000) (106,0.9300) (107,0.9100) (108,0.9400) (109,0.9300) (110,0.9300) (111,0.9500) (112,0.9800) (113,0.9200) (114,0.9500) (115,0.9700) (116,0.9600) (117,0.9800) (118,0.9600) (119,0.9900) (120,0.9900) (121,1.0000) (122,0.9600) (123,1.0000) (124,0.9900) (125,0.9800) (126,0.9800) (127,1.0000) (128,0.9800) (129,0.9900) (130,0.9700) (131,0.9800) (132,0.9800) (133,1.0000) (134,1.0000) (135,0.9800) (136,0.9900) (137,1.0000) (138,0.9800) (139,0.9900) (140,1.0000) (141,0.9900) (142,0.9900) (143,0.9900) (144,0.9900) (145,1.0000) (146,0.9900) (147,1.0000) (148,0.9900) (149,1.0000) (150,0.9800) (151,0.9800) (152,1.0000) (153,1.0000) (154,1.0000) (155,1.0000) (156,0.9900) (157,1.0000) (158,1.0000) (159,1.0000) (160,1.0000) 
};
\end{groupplot}
\end{tikzpicture}
\caption{Proportion of random digraphs on $n$ vertices with $m$ edges (self-loops allowed) satisfying the SCC-RH necessary condition (blue) and the peelable sufficient condition (red), for $n\in\{5,10,15,20\}$. For $n=5$ the proportions are calculated exactly and the proportion of stable graphs is also shown (green). For $n\in\{10,15,20\}$ the proportion is estimated over $100$ independent samples.}
\label{fig:gap-rh-peelable-large}
\end{figure*}

To see how the phase transition scales, we generate random digraphs of $n$ vertices with a prescribed number of edges $m$, sampled uniformly (without replacement) from all $n^2$ ordered pairs. For $n\in\{10,15,20\}$ and varying $m$, we generate $100$ digraphs and
record the proportions of samples that are SCC-RH and peelable (we are no longer able to assess stability due to its impractical computational cost), as shown in Figure~\ref{fig:gap-rh-peelable-large}. The smallest number of edges $m$ at which both proportions are above
$50\%$ is henceforth called the \emph{critical edge count}, and it equals $12,31,54,79$ for $n=5,10,15,20$. The maximum gap between SCC-RH and peelability occurs at edge count $m=10,27,47,67$ for $n=5,10,15,20$, and stays above $9\%$, rather than closing.
The results suggest that the gap between both properties does not diminish as $n$ grows.

\subsection{Running time}

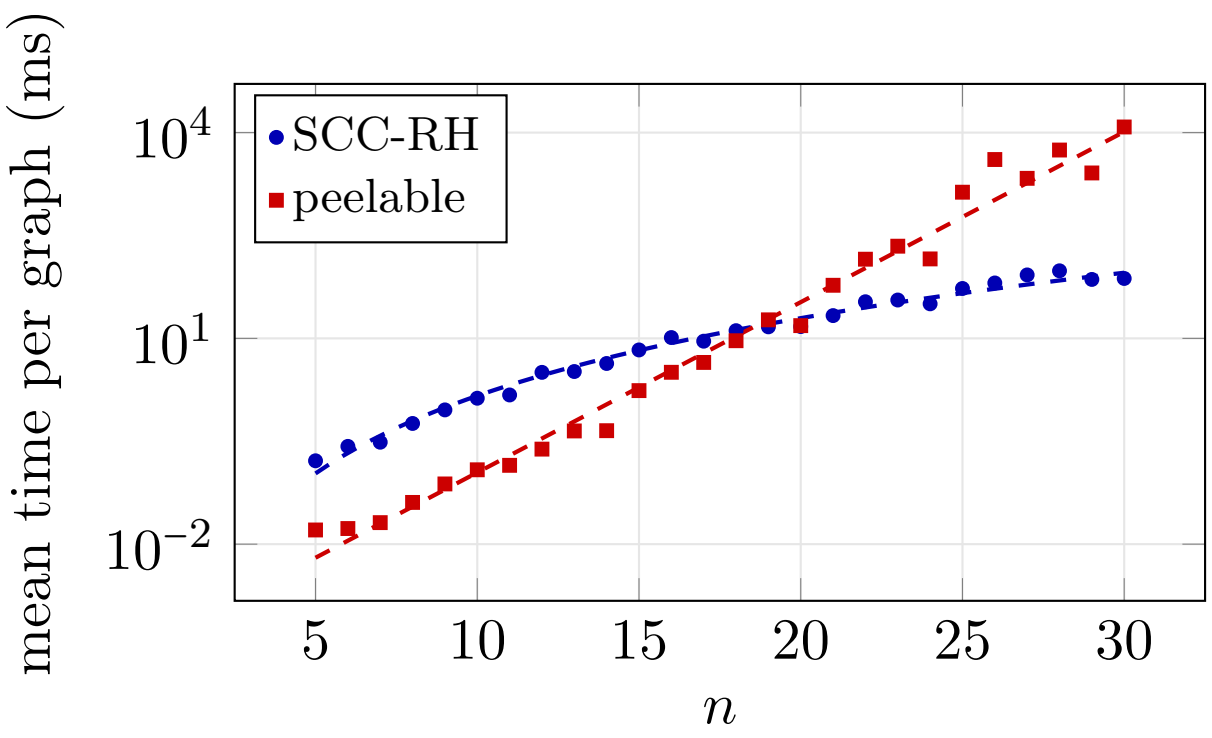
\begin{figure}[!ht]\centering
\begin{tikzpicture}
\begin{axis}[width=8cm, height=5cm, ymode=log,
  xlabel={$n$}, ylabel={mean time per graph (ms)},
  grid=both, grid style={gray!20},
  legend style={font=\footnotesize, at={(0.02,0.98)}, anchor=north west},
  legend cell align=left]
\addplot[mark=*, thick, color=blue!70!black, mark size=1pt, only marks] coordinates {(5,0.1642) (6,0.2664) (7,0.3064) (8,0.5750) (9,0.9067) (10,1.3413) (11,1.4966) (12,3.1970) (13,3.2881) (14,4.3060) (15,6.7841) (16,10.2853) (17,9.1021) (18,12.9542) (19,14.7382) (20,14.8767) (21,21.4905) (22,34.2146) (23,36.2511) (24,31.9137) (25,53.3988) (26,64.2545) (27,84.1495) (28,96.6399) (29,72.3587) (30,74.7370)};
\addplot[mark=square*, thick, color=red!80!black,mark size=1pt, only marks] coordinates {(5,0.0161) (6,0.0169) (7,0.0206) (8,0.0406) (9,0.0755) (10,0.1214) (11,0.1410) (12,0.2430) (13,0.4472) (14,0.4511) (15,1.7319) (16,3.2035) (17,4.4603) (18,9.2458) (19,18.6401) (20,15.3478) (21,59.4673) (22,142.7134) (23,220.8579) (24,144.0670) (25,1352.6530) (26,4052.1306) (27,2154.9526) (28,5571.5572) (29,2580.8032) (30,12078.4132)};

\addplot[domain=5:30, samples=60, dashed, thick, color=blue!70!black, forget plot]
  {0.000253*x^3.76}; 
\addplot[domain=5:30, samples=60, dashed, thick, color=red!80!black, forget plot]
  {0.000363*exp(0.572*x)};
  
\legend{SCC-RH, peelable}
\end{axis}
\end{tikzpicture}
\caption{Mean running time, per random digraph, of the two conditions, at the estimated critical edge count, as a function of $n$. Dashed lines are least-squares fits on the log scale: $\approx 2.5\cdot 10^{-4}\,n^{3.76}\ \mathrm{ms}$ for SCC-RH and $\approx 3.6\cdot 10^{-4}\,1.77^{n}\ \mathrm{ms}$ for peelable.}
\label{fig:timing-vs-n}
\end{figure}

We also measured the running time of the algorithms for SCC-RH and peelability, as the size of the graph varies. We estimate the critical edge count, via binary search on the number of edges $m$, as the smallest $m$ achieving at least 50\% SCC-RH graphs among $T=100$ samples. 
Figure~\ref{fig:timing-vs-n} reports the mean time (averaged over 100 samples) on random graphs at the estimated critical edge count. 
Least-squares fitting suggests a running time approximately proportional to $n^{3.76}$ for SCC-RH (Algorithm~\ref{alg:rh-property}) and $1.77^n$ for peelability (Algorithm~\ref{alg:peelability}). These results are consistent with our theoretical worst-case analysis. The cost of peelability overtakes that of the necessary condition at $n=19$ and, at $n=30$, it is about $160$ times slower.

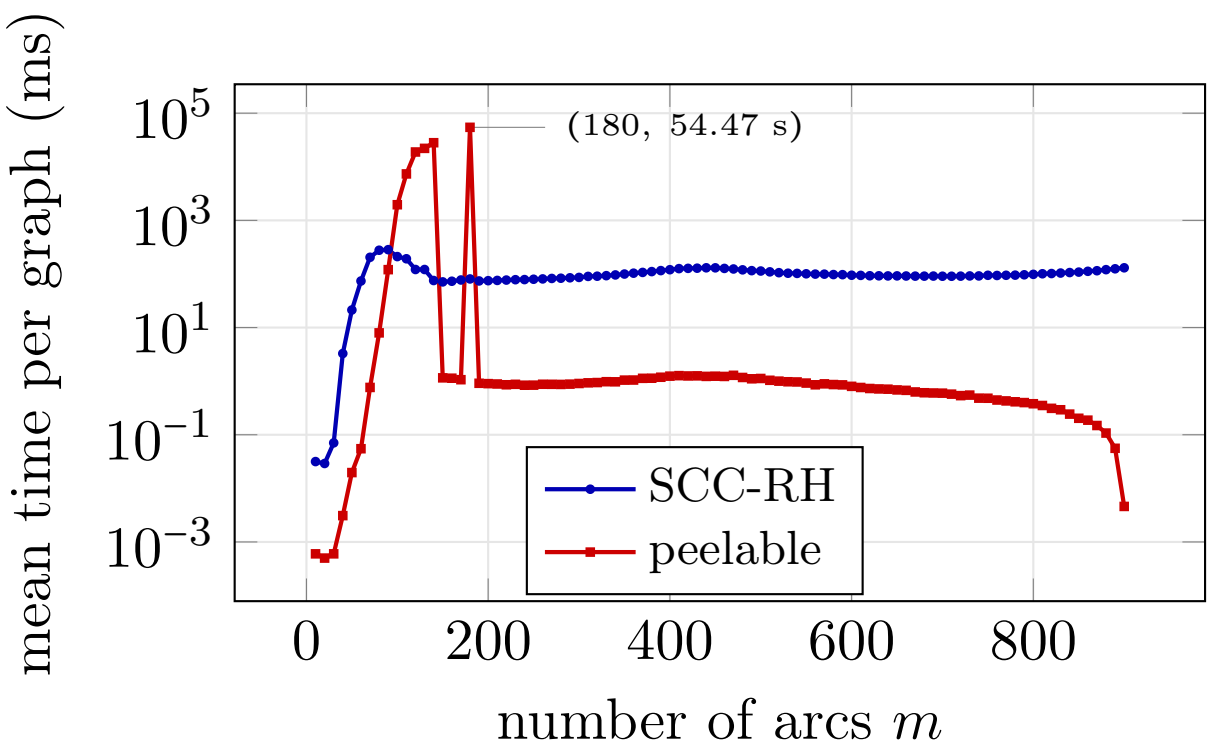
\begin{figure}[!ht]\centering
\begin{tikzpicture}
\begin{axis}[width=8cm, height=5cm, ymode=log,
  xlabel={number of arcs $m$}, ylabel={mean time per graph (ms)},
  grid=both, grid style={gray!20},
  legend style={font=\footnotesize, at={(0.3,0.3)}, anchor=north west},
  ytick={10^-3,10^-1,10^1,10^3,10^5},
  legend cell align=left]
\addplot[mark=*, thick, color=blue!70!black, mark size=.5pt] coordinates {(10,0.0316) (20,0.0291) (30,0.0706) (40,3.2824) (50,21.2969) (60,73.9542) (70,204.7653) (80,277.2789) (90,284.1594) (100,210.9336) (110,192.1427) (120,121.1889) (130,121.1879) (140,75.5055) (150,71.1778) (160,73.2021) (170,77.2077) (180,80.7678) (190,74.1101) (200,74.9688) (210,75.8348) (220,76.9945) (230,78.2139) (240,78.5779) (250,79.7430) (260,80.8827) (270,82.5156) (280,83.3101) (290,84.8502) (300,86.2172) (310,90.1018) (320,90.8755) (330,93.0033) (340,95.9444) (350,99.6914) (360,103.4778) (370,106.9006) (380,110.7752) (390,115.5460) (400,120.1098) (410,125.8297) (420,127.3120) (430,128.2025) (440,130.3610) (450,130.1652) (460,126.9739) (470,123.9291) (480,119.4239) (490,115.3656) (500,113.5198) (510,109.5520) (520,105.6935) (530,103.2730) (540,102.5199) (550,100.7464) (560,99.4487) (570,99.2225) (580,97.5233) (590,97.1281) (600,94.7516) (610,94.1224) (620,92.6353) (630,92.4175) (640,91.9718) (650,92.1351) (660,91.6682) (670,91.1581) (680,91.1065) (690,91.3903) (700,90.7792) (710,90.6772) (720,91.0771) (730,91.8984) (740,91.9641) (750,94.4910) (760,93.7540) (770,94.5511) (780,95.0952) (790,96.5307) (800,98.3814) (810,100.9063) (820,101.8671) (830,104.0376) (840,106.2841) (850,108.3702) (860,112.1487) (870,114.9652) (880,120.3529) (890,124.9147) (900,130.3552)};
\addplot[mark=square*, thick, color=red!80!black, mark size=.5pt] coordinates {(10,0.0006) (20,0.0005) (30,0.0006) (40,0.0031) (50,0.0197) (60,0.0549) (70,0.7639) (80,7.9652) (90,120.4280) (100,1959.1340) (110,7360.3058) (120,18877.1817) (130,22056.0365) (140,28119.6840) (150,1.1530) (160,1.1352) (170,1.0650) (180,54472.7962) (190,0.9131) (200,0.8989) (210,0.8857) (220,0.8567) (230,0.8736) (240,0.8403) (250,0.8426) (260,0.8756) (270,0.8739) (280,0.8693) (290,0.8804) (300,0.9049) (310,0.9306) (320,0.9422) (330,0.9768) (340,0.9725) (350,1.0458) (360,1.0494) (370,1.1296) (380,1.1349) (390,1.1834) (400,1.2401) (410,1.2674) (420,1.2507) (430,1.2561) (440,1.2277) (450,1.2379) (460,1.2248) (470,1.2857) (480,1.1719) (490,1.1069) (500,1.1222) (510,1.0383) (520,1.0007) (530,0.9748) (540,0.9662) (550,0.9215) (560,0.8556) (570,0.8911) (580,0.8645) (590,0.8524) (600,0.8004) (610,0.7625) (620,0.7297) (630,0.7160) (640,0.7087) (650,0.6865) (660,0.6741) (670,0.6291) (680,0.6080) (690,0.6018) (700,0.5973) (710,0.5726) (720,0.5369) (730,0.5494) (740,0.4820) (750,0.4800) (760,0.4447) (770,0.4280) (780,0.4121) (790,0.3984) (800,0.3791) (810,0.3517) (820,0.3112) (830,0.2927) (840,0.2428) (850,0.2030) (860,0.1872) (870,0.1488) (880,0.1076) (890,0.0559) (900,0.0046)};

\node[coordinate,pin=right :{\tiny $(180,\,\,54.47\ \mathrm{s})$}] at (axis cs:180,54472.7262){};

\legend{\footnotesize SCC-RH, \footnotesize peelable}
\end{axis}
\end{tikzpicture}
\caption{Mean running time, per random digraph, at $n=30$, as a function of the number of arcs $m$. The cost of peelability rises sharply around the critical edge count $m\approx 130$.}
\label{fig:timing-vs-m}
\end{figure}

Finally, we study how the number of edges influences the computation cost. 
For figure~\ref{fig:timing-vs-m}, we fix $n=30$, vary $m$, and report the mean time (averaged over 100 samples) on random digraphs with $m$ edges. 
While the cost of the necessary condition stays below~$1$ second throughout, the cost of peelability spikes to about $28$ seconds around $m=140$ (with an even larger outlier at $m=180$). 
We conclude that the necessary condition stays uniformly cheap across all densities, whereas peelability is easy on sufficiently sparse and sufficiently dense graphs but becomes drastically hard at intermediate densities.

\section{Conclusion}

We revisited the graphical stability conditions of \cite{belabbas2013sparse}
and showed that the bracket they place around stability is loose on both sides,
and that it does not tighten under strong connectivity. The sufficient
condition, a nested Hamiltonian decomposition, fails to be necessary already at
order three. The necessary condition, packability together with SCC-rootedness, fails to be sufficient already at order four. 
Refining the necessary condition from the
coefficients of the characteristic polynomial to the full Routh-Hurwitz table
detects the second example but remains strictly weaker than stability at order five. 

The two sides also differ in cost. 
Each of the necessary conditions reduces to
polynomial identity testing and lies in $\mathrm{RP}$, hence is not
$\mathrm{NP}$-hard unless $\mathrm{NP}=\mathrm{RP}$; in contrast, peelability
is $\mathrm{NP}$-complete. 
Our experiments trace this contrast in practice: the algorithm that checks the
necessary conditions is cheap at every size and density, while the algorithm that checks peelability is easy on sparse and on dense graphs but hard in an \mbox{intermediate-density} region.

These results leave open the exact complexity of deciding stability of a general sparse matrix space. We read the two-sided gap as evidence that a tractable criterion, should
one exist, will have to move beyond the cycle and connectivity structure that both of Belabbas' conditions rest on.

{\small
\bibliographystyle{plainnat}
\bibliography{bib.bib}
}

\vspace{4cm}

\end{document}